\documentclass[12pt,reqno]{amsart}

\usepackage{amssymb,amsmath, amsfonts, latexsym}
\usepackage{enumerate}
\newtheorem{thm}{Theorem}[section]
\newtheorem{cor}[thm]{Corollary}
\newtheorem{lem}[thm]{Lemma}
\newtheorem{prop}[thm]{Proposition}

\newtheorem{rem}{Remark}[section]

\newcommand{\dR}{\mathbb{R}}
\newcommand{\dP}{\mathbb{P}}

\newcommand{\superexp} {\quad\underset{\lambda_n}{\overset{\rm superexp}{\longrightarrow}} }

\newcommand{\superexpldp} {\quad\underset{n}{\overset{\rm superexp}{\longrightarrow}} }
\newcommand{\superexpmdp} {\quad \underset{b^2_n}{\overset{\rm superexp}{\longrightarrow}} }

\def\build#1_#2^#3{\mathrel{\mathop{\kern 0pt#1}\limits_{#2}^{#3}}}

 \numberwithin{equation}{section}

\begin{document}

\title[LDP of the realized (co-)volatility]{Large deviations of the realized (co-)volatility vector}

\author{Hac\`ene Djellout}
\email{Hacene.Djellout@math.univ-bpclermont.fr}
\address{Laboratoire de Math\'ematiques, CNRS UMR 6620, Universit\'e Blaise Pascal, Avenue des Landais, 63177 Aubi\`ere, France.}

\author{Arnaud Guillin}
\email{Arnaud.Guillin@math.univ-bpclermont.fr}
\address{Laboratoire de Math\'ematiques, CNRS UMR 6620, Universit\'e Blaise Pascal, Avenue des Landais, 63177 Aubi\`ere, France.}

\author{Yacouba Samoura}
\email{Yacouba.Samoura@math.univ-bpclermont.fr}
\address{Laboratoire de Math\'ematiques, CNRS UMR 6620, Universit\'e Blaise Pascal, Avenue des Landais, 63177 Aubi\`ere, France.}

\keywords{Realised Volatility and covolatility, large  deviations, diffusion, discrete-time observation}

\date{\today}
\begin{abstract} Realized statistics based on high frequency returns have become very popular in financial economics. In recent years, different non-parametric estimators of the variation of a log-price process have appeared. These were developed by many authors and were motivated by the existence of complete records of price data. Among them are the realized quadratic (co-)variation which is perhaps the most well known example, providing a consistent estimator of the integrated (co-)volatility  when the logarithmic price process is continuous. Limit results such as the weak law of large numbers or the central limit theorem have been proved in different contexts. In this paper, we propose to study the large deviation properties of realized (co-)volatility (i.e., when the number of high frequency observations in a fixed time interval increases to infinity. More specifically, we consider a bivariate model with synchronous observation schemes and correlated Brownian motions of the following form:  $dX_{\ell,t} = \sigma_{\ell,t}dB_{\ell,t}+b_{\ell}(t,\omega)dt$ for $\ell=1,2$, where $X_{\ell}$ denotes the log-price, we are concerned  with  the  large deviation estimation of the vector 
$V_t^n(X)=\left(Q_{1,t}^n(X), Q_{2,t}^n(X), C_{t}^n(X)\right)$ where $Q_{\ell,t}^n(X)$ and $C_{t}^n(X)$ represente the estimator of the quadratic variational processes $Q_{\ell,t}=\int_0^t\sigma_{\ell,s}^2ds$ and the integrated covariance $C_t=\int_0^t\sigma_{1,s}\sigma_{2,s}\rho_sds$ respectively, with $\rho_t=cov(B_{1,t}, B_{2,t})$. Our main motivation is to improve upon the existing limit theorems. 
Our large deviations results can be used to evaluate and approximate tail probabilities of realized (co-)volatility. As an application we provide the large deviation for the standard dependence measures between the two assets returns such as the realized regression coefficients up to time $t$, or the realized correlation. Our study should contribute to the recent trend of research on the (co-)variance estimation problems, which are quite often discussed in high-frequency financial data analysis. 
\end{abstract}

\maketitle

\vspace{-0.5cm}

\begin{center}
\textit{AMS 2000 subject classifications: 60F10, 60G42, 62M10, 62G05.}
\end{center}

\medskip

\maketitle


\section{Introduction, Model and Notations}
In the last decade there has been a considerable development of the asymptotic theory for processes observed at a high frequency. This was mainly motivated by financial applications, where the data, such as stock prices or currencies, are observed very frequently. 

Asset returns covariance and its related statistics play a prominent role in many important theoretical as well as practical problems in finance. Analogous to the realized volatility approach, the idea of employing high frequency data in the computation of daily (or lower frequency) covariance between two assets leads to the concept of realized covariance (or covariation). The key role of quantifying integrated (co-)volatilities in portfolio optimization and risk management has stimulated an increasing interest in estimation methods for these models.

\vspace{10pt}
It is quite natural to use the asymptotic framework when the number of high frequency observations in a fixed time interval (say, a day) increases to infinity. Thus Barndorff-Nielsen and Shephard \cite{BNS} established a law of large numbers and  the corresponding fluctuations for realized volatility, also extended to more general setups and statistics by Barndorff-Nielsen et al. \cite{BNGJS} and \cite{BNGJPS}. Dovonon, Gon{\c{c}}alves, and Meddahi \cite{DGM} considered Edgeworth expansions for the realized volatility statistic and its bootstrap analog. These results are crucial to explore asymptotic behaviors of realized (co-)volatility, in particular around the center of its distribution. There are also different estimation approaches for the integrated covolatility in multidimensional models and limit theorem, and we can refer to Barndorff-Nielsen et al. \cite{BNSW} and \cite{BNGJS} where the authors present, in an unified way, a weak law of large numbers and a central limit theorem for a general estimator, called realized generalized bipower variation. 

\vspace{5pt}
For related work concerning bivariate case under a non-synchronous sampling scheme, see Hayashi and Yoshida \cite{Hayashi1}, Bibinger \cite{B},  Dalalyan and Yoshida \cite{Dalalyan1}, see also A\"\i t-Sahalia et al. \cite{Sahalia1} and the references therein. Estimation of the covariance of log-price processes in the presence of market microstructure noise, we refer to Bibinger and Rei{\ss} \cite{BR}, Robert and Rosenbaum \cite{RM}, Zhang et al. \cite{ZMA2} and \cite{ZMA1}.  See also Gloter, or Comte et al. \cite{CGCR} for non parametric estimation in the case of a stochastic volatility model.

\vspace{10pt}
We model the evolution of an observable state variable by a stochastic process $X_t=(X_{1,t}, X_{2,t}), t\in [0, 1].$ In financial applications, $X_t$ can be thought of as the short interest rate, a foreign exchange rate, or the logarithm of an asset price or of a stock index. Suppose both $X_{1,t} $ and $X_{2,t} $ are defined on a filtered probability space $(\Omega, \mathcal{F}, (\mathcal{F}_t), \mathbb{P})$ and  follow an It\^o process, namely,
\begin{equation}
\label{AR}
\vspace{1ex}
\left\{
\begin{array}[c]{ccccc}
dX_{1,t} & = & \sigma_{1,t}dB_{1,t} & + & b_{1}(t,\omega)dt\\
dX_{2,t} & = & \sigma_{2,t}dB_{2,t} & + & b_{2}(t,\omega)dt
\end{array}
\right.
\end{equation}
where $B_1$ and $B_2$ are standart Brownian motions, with correlation $Corr(B_{1,t},B_{2,t})=\rho_t$. We can write $dB_{2,t}=\rho_tdB_{1,t}+\sqrt{1-\rho_t^2}dB_{3,t},$ where $B_1=(B_{1,t})_{t\in [0,1]}$ and  $B_3=(B_{3,t})_{t\in [0,1]}$ are independent Brownian processes. 
\vspace{10pt}

We will suppose of course existence and uniqueness of strong solutions, and in what follows, the drift coefficient $b_1$ and $b_2$ are assumed to satisfy an uniform linear growth condition and we limit our attention to the case when $\sigma_1$, $\sigma_2$ and $\rho$ are deterministic functions. The functions $\sigma_{\ell}$, $\ell=1,2$ take positive values while $\rho$ takes values in the interval $]-1, 1[.$ 

\vspace{10pt}

In this paper, our interest is to estimate the (co-)variation vector 
\begin{equation}\label{vector}
[V]_t=([X_1]_t,[X_2]_t, \left\langle X_1,X_2 \right\rangle_{t})^T
\end{equation}
between two returns in a fixed time period $[0; 1]$ when $X_{1,t}$ and $X_{2,t}$ are observed synchronously, $[X_{\ell}]_t,\ell=1,2$ represente
the quadratic variational process of $X_{\ell}$ and $\left\langle X_1,X_2 \right\rangle_{t}$ the (deterministic) covariance of $X_1$ and $X_2$:
$$[X_{\ell}]_t=\int_0^t \sigma_{\ell,s}^2 \mathrm ds, \qquad \left\langle X_1,X_2 \right\rangle_{t}=\int_0^t \sigma_{1,s}\sigma_{2,s}\rho_s \mathrm ds.$$

Inference for (\ref{vector}) is a well-understood problem if $X_{1,t}$ and $X_{2,t}$ are observed simultaneously. Note that $X_{1,t}$ and $X_{2,t}$ are not observed in continuous time but we have only discrete time observations. Given discrete equally space observation
$(X_{1,t_k^n},X_{2,t_k^n},k=1,\cdots,n)$ in the interval $[0,1]$ (with $\,\,t_k^n = k/n)$, a limit theorem in stochastic processes states that  
$$V_t^n(X)=\left(Q_{1,t}^n(X), Q_{2,t}^n(X), C_{t}^n(X)\right)^T$$
commonly called realized (co-)variance, is a consistent estimator for $[V]_t$, with, for $\ell=1,2$
$$Q_{\ell,t}^n(X)= \sum_{k=1}^{[nt]}\left(\Delta_{k}^nX_{\ell}\right)^2\qquad C_{t}^n(X) = \sum_{k=1}^{[nt]}\left( \Delta_{k}^nX_{1} \right)\left( \Delta_{k}^nX_{2}\right)$$
where $[x]$ denote the integer part of $x \in \mathbb{R}$ and $\Delta_{k}^nX_{\ell}= X_{\ell,t_{k}^n}-X_{\ell,t_{k-1}^n}$.

\vspace{10pt}

When the drift $b_{\ell}(t,\omega)$ is known, we can consider the following variant 
$$V_t^n(X-Y)=\left(Q_{1,t}^n(X-Y) , Q_{2,t}^n(X-Y) , C_{t}^n(X-Y)\right)^T$$ 
with for $\ell=1,2$ and $Y_{\ell,t}:= \int_{0}^{t}b_{\ell}(t,\omega) \mathrm dt,$
$$Q_{\ell,t}^n(X-Y) = \sum_{k=1}^{[nt]}\left(\Delta_{k}^nX_{\ell}-\Delta_{k}^nY_{\ell}\right)^2,$$
and
$$C_{t}^n(X-Y) = \sum_{k=1}^{[nt]}\left(\Delta_{k}^nX_{1}-\Delta_{k}^nY_{1}\right)\left(\Delta_{k}^nX_{2}-\Delta_{k}^nY_{2}\right).$$
\vspace{10pt}

If the drift $b_{\ell}(t,\omega):=b_{\ell}(t,X_{1,t}(\omega),X_{2,t}(\omega))$, where $b_{\ell}(t,x_1,x_2)$ is some deterministic function (a current situation), $X_t=(X_{1,t},X_{2,t})$  verifies
\begin{equation}
\label{AR22}
\vspace{1ex}
\left\{
\begin{array}[c]{ccccc}
dX_{1,t} & = & \sigma_{1,t}dB_{1,t} & + & b_{1}(t,X_t)dt\\
dX_{2,t} & = & \sigma_{2,t}dB_{2,t} & + & b_{2}(t,X_t)dt.
\end{array}
\right.
\end{equation}

When $b_{\ell}(t,x)$ is known, and only the sample $X_{t_{k-1}^n}=(X_{1,t_{k-1}^n},X_{2,t_{k-1}^n}),k=0 \cdots n$ is observed, we can also consider the following estimator $\tilde{V_t^n}(X)=\left(\tilde{Q_{1,t}^n}(X), \tilde{Q_{2,t}^n}(X), \tilde{C_{t}^n}(X)\right)^T$ with for $\ell=1,2$
$$\tilde{Q_{\ell,t}^n}(X) = \sum_{k=1}^{[nt]}\left(\Delta_k^nX_{\ell}-b_{\ell,t_{k-1}^n}(X_{t_{k-1}^n})(t_{k}^n-t_{k-1}^n)\right)^2,$$
$$\tilde{C_{t}^n} (X)= \sum_{k=1}^{[nt]}\left(\Delta_k^nX_1-b_{1,t_{k-1}^n}(X_{t_{k-1}^n})(t_{k}^n-t_{k-1}^n)\right)\left(\Delta_k^nX_2-b_{2,t_{k-1}^n}(X_{t_{k-1}^n})(t_{k}^n-t_{k-1}^n)\right).$$
\vspace{10pt}

In the aforementionned papers, and under quite weak assumptions, it is proved the following consistency
$$V^n_1(X),\, V^n_1(X-Y),\, \tilde{V^n_1}(X)\longrightarrow [V]_1\qquad a.s.$$
and the corresponding fluctuations
$$\sqrt{n}(V^n_1(X)-[V]_1), \,\sqrt{n}(V^n_1(X-T)-[V]_1), \,\sqrt{n}(\tilde{V^n_1}(X)-[V]_1), \stackrel{\mathcal L}{\longrightarrow}{\mathcal N}(0,\Sigma).$$

The purpose of this paper is to furnish some further trajectorial estimations about the estimator $V_{.}^n$, deepening the law of large numbers and central limit theorem. More
precisely, we are interested in the estimation of 
$$\mathbb{P}\left(\dfrac{\sqrt{n}}{b_n}(V_{.}^n(X)-[V]_{.}) \in A\right),$$
where A is a given domain of deviation, and $(b_n)_{n\geqslant0}$ is some sequence denoting the scale of the deviation. 

When $b_n=1,$ this is exactly the estimation of the central limit theorem. When $b_n=\sqrt{n},$ it becomes the large deviations. And when $1 << b_n << \sqrt{n},$ it is called moderate deviations. In other words, the moderate deviations investigate the convergence speed between the large deviations and central limit theorem.
\vspace{10pt}

The large deviations and moderate deviations problems arise in the theory of statistical inference quite naturally. For estimation of unknown parameters and functions, it is first of all important to minimize the risk of wrong decisions implied by deviations of the observed values of estimators from the true values of parameters or functions to be estimated. Such important errors are precisely the subject of large deviation theory. The large deviation and moderate deviation results of estimators can provide us with the rates of convergence and an useful method for constructing asymptotic confidence intervals.
\vspace{10pt}

The aim of this paper is then to focus on the large and moderate deviation estimations of the estimators of volatility and co-volatility.  Despite the fact that these statistics are nearly 20 years old, there has been remarkably few result in this direction, it is a surprise to us. The answer may however be the following: the usual techniques (such as G\"artner-Ellis method) do not work and a very particular treatment has to be considered for this problem. Recently, however, some papers considered the unidimensional case. Djellout et al. \cite{Djellout1} and recently Shin and Otsu \cite{OS} obtained the large and moderate deviations for the realized volatility.  In the bivariate case Djellout and Yacouba \cite{Djellout2}, obtained the large and moderate deviations for the realized covolatility. The large deviation for threshold estimator for the constant volatility was established  by Mancini \cite{Mancini1} in jumps case. And the moderate deviation for threshold estimator for the quadratic variational process was derived by Jiang \cite{Jiang1}. Let us mention that the problem of the large deviation for threshold estimator vector, in the presence of jumps,  will be considered in a forthcoming paper,  consistency, efficience and robustnesse were proved in Mancini and Gobbi \cite{Mancini2}. The case of asynchronous sampling scheme, or in the presence of micro-structure noise is also outside the scope of the present paper but are currently under investigations.

\vspace{10pt}
Two economically interesting functions of the realized covariance vector are the realized correlation and the realized regression coefficients. In particular, realized regression coefficients are obtained by regressing high frequency returns for one asset on high frequency returns for another asset. When one of the assets is the market portfolio, the result is a realized beta coefficient. A beta coefficient measures the assets systematic risk as assessed by its correlation with the market portfolio. Recent examples of papers that have obtained empirical estimates of realized betas include Andersen, Bollerslev, Diebold and Wu \cite{ABDW}, Todorov and Bollerslev \cite{TB}, Dovonon, Gon{\c{c}}alves and Meddahi \cite{DGM}, Mancini and Gobbi \cite{Mancini2}.

\vspace{5pt}
Let us stress that large deviations for the realized correlation can not be deduced from unidimensional quantities and were thus largely ignored. As an application of our main results, we provide a large and moderate deviation principle for the realized correlation and the realized regression coefficients in some special cases. The realized regression coefficient from regressing is $\beta_{\ell,t}^n(X)=\frac{C_t^n(X)}{Q^n_{\ell,t}(X)}$ which consistently estimates $\beta_{\ell,t}=\frac{C_t}{Q_{\ell,t}}$ and the realized correlation coefficient is $\varrho_t^n(X)=\frac{C_t^n(X)}{\sqrt{Q^n_{1,t}(X)Q^n_{2,t}(X)}}$ which estimates $\varrho_t=\frac{C_t}{\sqrt{Q_{1,t}Q_{2,t}}}.$  The application will be based essentially on an application of the delta method, developped by Gao and Zhao (\cite{Gao1}).

\vspace{10pt}
As in Djellout et al. \cite{Djellout1}, Shin and Otsu \cite{OS}, it should be noted that the proof strategy of G\"artner and Ellis large deviation theorem can not be adapted here int he large deviations case. We will encounter the same technical difficulties as in the papers of Bercu et al. \cite{BGR} and Bryc and Dembo \cite{BD} where they established the large deviation principle for quadratic forms of Gaussian processes. Since we cannot determine the limiting behavior of the cumulant generating function at some boundary point, we will use an other approach based on the results of Najim \cite{Najim1}, \cite{Najim2} and \cite{Najim3}, where the steepness assumption concerning the cumulant generating function is relaxed. It has to be noted that the form of the large deviations rate function is also original: at the process level, and because of the weak exponential integrability of $V_t^n$, a correction (or extra) term appears in rate function, a phenomenon first discovered by Lynch and Sethuraman \cite{LS}. 
 
\vspace{10pt}
To be complete, let us now recall some basic definitions of the large deviations theory (c.f \cite{Demzei}). Let $(\lambda_n)_{n\ge 1}$ be a sequence of nonnegative real number such that $\lim_{n \to \infty}\lambda_n=+\infty$.   We say that a sequence of a random variables  $(M_{n})_{n}$ with topological state space $(S,  \mathbb{S})$, where $\mathbb{S}$ is a $\sigma-algebra$ on $S$, satisfies a large deviation principle with speed $\lambda_{n}$ and rate function $I : S \rightarrow [0,+\infty]$ if, for each $A \in  \mathbb{S}$,
\begin{equation*}
-\inf_{x \in A^{o}} I(x) \leq \liminf_{n \rightarrow \infty} \frac{1}{\lambda_n} \log \dP \Big( M_{n} \in A \Big) \leq \limsup_{n \rightarrow \infty} \frac{1}{\lambda_n} \log \dP \Big( M_{n} \in A \Big) \leq -\inf_{x \in \bar{A}} I(x)
\end{equation*}
where $A^{o}$ and $\bar{A}$ denote the interior and the closure of $A$, respectively. 

The rate function $I$ is lower semicontinuous, \textit{i.e.} all the sub-level sets $\{ x \in S ~ \vert ~ I(x) \leq c \}$ are closed, for $c \geq 0$. If these level sets are compact, then $I$ is said to be a \textnormal{good rate function}. When the speed of the large deviation principle correspond to the regime between the central limit theorem and the law of large numbers, we talk of moderate deviation principle.

\vspace{10pt}

\noindent{\bf Notations.} \textit{In the whole paper, for any matrix $M$, $M^T$ and $\Vert M \Vert$ stand for the transpose and the euclidean norm of $M$, respectively. For any square matrix $M$, $\det(M)$ is the determinant of $M$. Moreover, we will shorten \textnormal{large deviation principle} by LDP and \textnormal{moderate deviation principle} by MDP. 
We denote by $\langle\cdot,\cdot \rangle$ the usual scalar product. For any process $Z_t$, $\Delta_k^nZ$ stands for the increment $Z_{t_k^n}-Z_{t_{k-1}^n}$.
In addition, for a sequence of random variables $(Z_{n})_{n}$ on $\dR^{d \times p}$, we say that $(Z_{n})_{n}$ converges $(\lambda_{n})-$superexponentially fast in probability to some random variable $Z$ if, for all $\delta > 0$,
\begin{equation*}
\limsup_{n \rightarrow \infty} \frac{1}{\lambda_n} \log \dP\Big( \left\Vert Z_{n} - Z \right\Vert > \delta \Big) = -\infty.
\end{equation*}
This \textnormal{exponential convergence} with speed $\lambda_{n}$ will be shortened as
\begin{equation*}
Z_{n} \superexp Z.
\end{equation*}}

\vspace{10pt}
The article is arranged in three upcoming sections and an appendix comprising some theorems used intensively in the paper, we have included them here for completeness. Section 2 is devoted to our main results on the LDP and MDP for the (co-)volatility vector. In Section 3, we  deduce applications for the realized correlation and the realized regression coefficients, when $\sigma_{\ell}$, for $\ell=1,2$ are constants. In section 4, we give the proof of these theorems.

\section{Main results}

Let $X_t=(X_{1,t},X_{2,t})$ be given by \eqref{AR}, and $Y_t=(Y_{1,t},Y_{2,t})$ where for $\ell =1,2\quad$$Y_{\ell,t}:=\int_0^tb_{\ell}(t,\omega)dt$. We introduce the following conditions 

\vspace{10pt}
{\bf (B)} for $\ell=1,2$ $b(\cdot,\cdot)\in L^\infty(dt\otimes \mathbb{P})$

\vspace{10pt}
{\bf (LDP)} Assume that for $\ell=1,2$ 
\begin{itemize}
\item $\sigma_{\ell,t}^2(1-\rho_t^2)$ and $\sigma_{1,t}\sigma_{2,t}(1-\rho_t^2)$ $\in$ $L^\infty([0,1],dt)$.
\item the functions $t \to \sigma_{\ell,t}$ and $t \to \rho_t$ are continuous.
\end{itemize}

\vspace{10pt}
{\bf (MDP)} Assume that for $\ell=1,2$ 
\begin{itemize}
\item $\sigma_{\ell,t}^2(1-\rho_t^2)$ and $\sigma_{1,t}\sigma_{2,t}(1-\rho_t^2)$ $\in$ $L^2([0,1],dt)$.
\item Let $(b_n)_{n\geqslant1}$ be a sequence of positive
numbers such that
\begin{eqnarray}\label{Max_condition}
&&b_n \xrightarrow[{n \to \infty}]\quad  \infty \quad {\rm and} \quad \dfrac{b_n}{\sqrt{n}} \xrightarrow[{n \to \infty}]\quad 0\nonumber \\
&&{\rm and \,\,for}\quad \ell=1,2\quad\sqrt{n}b_n\max_{1 \leqslant k \leqslant n}\int_{(k-1)/n}^{k/n} \sigma_{\ell,t}^2 \mathrm dt \xrightarrow[{n \to \infty}]\quad 0.
\end{eqnarray}
\end{itemize}
\vspace{10pt}

We introduce the following function, which will play a crucial role in the calculation of the moment generating function: for $-1<c<1$ let for any $\lambda=(\lambda_1,\lambda_2,\lambda_3)\in \mathbb{R}^3$ 
\begin{equation}\label{ASSA}
P_c(\lambda):=
  \left\{
  \begin{array}{lcl}
  \vspace{0.3cm}
 -\dfrac{1}{2} \log\left(\dfrac{(1-2\lambda_1(1-c^2))(1-2\lambda_2(1-c^2))-(\lambda_3(1-c^2)+c)^2}{1-c^2}\right)\\
 \vspace{0.3cm}
\qquad \qquad\qquad \qquad \qquad\qquad if\qquad\lambda\in {\mathcal D}\\
 
  +\infty, \quad otherwise
  \end{array}
  \right.
\end{equation}
where 
\begin{equation}\label{defD}
\displaystyle {\mathcal D}_c=\left\{\lambda\in \mathbb{R}^3,\,\,\max_{\ell=1,2}\lambda_{\ell} < \dfrac{1}{2(1-c^2)}\,\,{\rm and}\,\prod_{\ell =1}^2\left(1-2\lambda_{\ell}(1-c^2)\right)> \left(\lambda_3(1-c^2)+c\right)^2\right\}.
\end{equation}
\vspace{10pt}

Let us present now the main results.

\subsection{Large deviation}

Our first result is about the large deviation of $V_{1}^n (X)$, i.e. at fixed time.
\begin{thm}\label{LDP}Let $t=1$ be fixed.
\begin{itemize}
 \item[(1)] For every $\lambda=(\lambda_1,\lambda_2,\lambda_3) \in \mathbb{R}^3,$
 \begin{equation*}
  \lim_{n \to \infty}\dfrac{1}{n}\log\mathbb{E}(\exp(n\left\langle \lambda,V_1^n(X-Y) \right\rangle))=\Lambda(\lambda):=\int_0^1 P_{\rho_t}(\lambda_1\sigma_{1,t}^2,\lambda_2\sigma_{2,t}^2,\lambda_3\sigma_{1,t}\sigma_{2,t}) \mathrm dt,
 \end{equation*}

where the function $P_c$ is given in (\ref{ASSA}).
\vspace{10pt}
 
 \item[(2)]  Under the conditions {\bf (LDP)} and {\bf (B)}  , the sequence $V_{1}^n(X)$ satisfies the LDP on $\mathbb{R}^3$ with speed $n$ and with the good rate function given
  by the legendre transformation of $\Lambda$, that is
  \begin{equation}\label{LDPtaux1}
  I_{ldp}(x)=\sup_{\lambda\in\mathbb{R}^3}\left(\left\langle \lambda,x \right\rangle - \Lambda(\lambda)\right).\end{equation}
\end{itemize}
\end{thm}

Let us consider the case where diffusion and correlation coefficients are constant, the rate function being easier to read (see also \cite{Shen1} in the purely Gaussian case, i.e. $b=0$). Before that let us introduce the function $P^*_c$ which is the Legendre transformation of $P_c$ given in (\ref{ASSA}), for all $x=(x_1,x_2,x_3)$ 
\begin{equation}\label{Legendre}
P^*_c(x):=
  \left\{
  \begin{array}{lcl}
  \vspace{0.3cm}
  \log\left(\dfrac{\sqrt{1-c^2}}{\sqrt{x_1x_2-x_3^2}}\right)-1+\dfrac{x_1+x_2-2cx_3}{2(1-c^2)}\\
 \vspace{0.3cm}
\qquad \qquad\qquad \qquad if\qquad x_1>0,\,\,x_2>0,\,\,x_1x_2>x_3^2\\
 
  +\infty, \quad otherwise.
  \end{array}
  \right.
\end{equation}
\begin{cor} \label{LDPconstant}
We assume that  for $\ell=1,2$ $\sigma_{\ell}$ and $\rho$ are constants. Under the condition {\bf (B)},  we obtain that $V_{1}^n(X)$ satisfies the LDP on $\mathbb{R}^3$ with speed $n$ and with the good rate function $I_{ldp}^{V}$ given by
\begin{equation}\label{LDPtauxconstantV}
I_{ldp}^{V}(x_1,x_2,x_3)=P^*_{\rho}\left(\frac{x_1}{\sigma_1^2},\frac{x_2}{\sigma_2^2}, \frac{x_3}{\sigma_1\sigma_2 }\right),
\end{equation}
where $P^*_c$ is given in (\ref{Legendre}).
\end{cor}
\vspace{10pt}

Now, we shall extend the Theorem \ref{LDP} to the process-level large deviations, i.e. for trajectories $(V^n_t(X))_{0\le t\le1}$, which is interesting 
from the viewpoint of non-parametric statistics.
\vspace{5pt}

Let $\mathcal BV([0,1],\mathbb{R}^3)$ (shortened in $\mathcal BV$) be the space of functions of bounded variation on $[0,1]$. We identify $\mathcal BV$ with $\mathcal M_3([0,1])$, the set of vector measures with value in $\mathbb{R}^3$. This is done in the usual manner:  to $f\in \mathcal BV$ there corresponds $\mu^f$ caracterized by  $\mu^f([0,t])=f(t)$. Up to this identification, $\mathcal C_3([0,1])$ the set of $\mathbb{R}^3$-valued continuous bounded functions on $[0,1]$), is the topological dual of $\mathcal BV$. We endow $\mathcal BV$ with the weak-* convergence topology $\sigma(\mathcal BV,\mathcal C_3([0,1]))$ (shortened $\sigma_w$) and with the associated Borel $\sigma-$field ${\mathcal B}_w$. Let $f\in \mathcal BV$ and $\mu^f$ the associated measure in $\mathcal M_3([0,1])$. Consider the Lebesgue  decomposition of $\mu^f$, $\mu^f=\mu_a^f+\mu_s^f$ where $\mu_a^f$ denotes the absolutely continous part of  $\mu^f$ with respect to $dx$ and $\mu_s^f$ its singular part. We denote by $f_a(t)=\mu_a^f([0,t])$ and by $f_s(t)=\mu_s^f([0,t])$.
\vspace{10pt}

\begin{thm}
 \label{LDP2} 
Under the conditions {\bf (LDP)} and {\bf (B)}, the sequence $V_{\cdot}^n(X)$ satisfies the LDP on $\mathcal BV$ with speed $n$ and with the rate function $J_{ldp}$ given for any $f=(f_1,f_2,f_3)\in \mathcal BV$ by
 \begin{eqnarray*}
 J_{ldp}(f)&=&\int_0^1P^*_{\rho_t}\left(\frac{f_{1,a}^{'}(t)}{\sigma_{1,t}^2},\frac{f_{2,a}^{'}(t)}{\sigma_{2,t}^2}, \frac{f_{3,a}^{'}(t)}{\sigma_{1,t}\sigma_{2,t}}\right) \mathrm dt\\ 
\qquad&&+ \int_0^1 \dfrac{\sigma_{2,t}^2f_{1,s}^{'}(t)+\sigma_{1,t}^2f_{2,s}^{'}(t)-2\rho_t\sigma_{1,t}\sigma_{2,t}f_{3,s}^{'}(t)}{2\sigma^2_{1,t}\sigma^2_{2,t}(1-\rho^2_t)}1_{[t;f_{1,s}^{'}>0,f_{2,s}^{'}>0,(f_{3,s}^{'})^2<f_{1,s}^{'}f_{2,s}^{'}]} \,\,\mathrm d\theta(t),
 \end{eqnarray*}
 where $P^*_c$ is given in (\ref{Legendre})  and $\theta$ is any real-valued nonnegative measure with respect to which $\mu_s^f$ is absolutely continuous and $f_s'=d\mu_s^f/d\theta=(f_{1,s}^{'},f_{2,s}^{'},f_{3,s}^{'})$.
\end{thm}
\vspace{10pt}

\begin{rem} Note that the definition of $f'_s$ is $\theta-$dependent. However, by homogeneity, $J_{ldp}$ does not depend upon $\theta$. One can choose $\theta=|f_{1,s}|+|f_{2,s}|+|f_{3,s}|$, with $|f_{i,s}|=f_{i,s}^++f_{i,s}^-$, where $f_{i,s}=f_{i,s}^+-f_{i,s}^-$ by the Hahn-Jordan decomposition.
\end{rem}

\begin{rem} As stated above, the problem of the LDP for $Q_{\ell,\cdot}^n(X)$ and  $C_{\cdot}^n(X)$ was alreay studied by Djellout et al. \cite{Djellout1} and \cite{Djellout2}, and the rate function is given explicitly in the last case. This is the first time that the LDP is investigated for the vector of the (co-)volatility.
\end{rem}

\begin{rem}
By using the contraction principle, and if $\sigma_{\ell}$ is strictly positive, we may find back the result of \cite{Djellout1}, i.e. that $Q^n_{\ell,\cdot}$ satisfies a LDP with speed $n$ and rate function
$$J_{ldp}^{\sigma_{\ell}}(f)=\int_0^1{\mathcal P}^*\left(\frac{f'_a(t)}{\sigma^2_{\ell,t}}\right)dt+\frac12\int_0^1\frac1{\sigma^2_{\ell,t}}d|f_s|(t)$$
where ${\mathcal P}^*(x)=\frac12(x-1-\log(x))$ when $x$ is positive and infinite if non positive, using the same notation as in the theorem (with $\theta=|f_s|$). One may also obtain the LDP for $C^n_\cdot$ by the contraction principle, recovering the result of Djellout-Yacouba \cite{Djellout2} (see there for the quite explicit complicated rate function). 
\end{rem}
\begin{rem}
The continuity assumptions in {\bf (LDP)} of $\sigma_{\ell,\cdot}$ and $\rho_\cdot$ is not necessary, but in this case we have to consider another strategy of the proof, more technical and relying on Dawson-G\"artner type theorem, which moreover does not enable to get other precision on the rate function that the fact it is a good rate function.\\ However it is not hard to adapt our proof to the case where $\sigma_{\ell,\cdot}$ and $\rho_\cdot$ have only a finite number of discontinuity points (of the first type). This can be done by applying the previous theorem to each subinterval where all functions are continuous and using the independence of the increments of $V^n_t(X-Y)$.
\end{rem}
\vspace{10pt}


\subsection{Moderate deviation}

Let us now considered the intermediate scale between the central limit theorem and the law of large numbers.
\begin{thm}\label{MDP1} For t=1 fixed. Under the conditions {\bf (MDP)} and {\bf (B)}  , the sequence
 $$\dfrac{\sqrt{n}}{b_n}\left(V_{1}^n(X)-[V]_{1}\right)$$ 
 satisfies the LDP on $\mathbb{R}^3$ with speed $b_n^2$ and with the rate function given by
 \begin{equation}\label{MDPtaux1}
 I_{mdp}(x)=\sup_{\lambda\in\mathbb{R}^3}\left(\left\langle \lambda,x \right\rangle - \dfrac{1}{2}\left\langle\lambda, \Sigma_1\cdot\lambda\right\rangle\right)=\frac{1}{2}\left\langle x, \Sigma_1^{-1}\cdot x\right\rangle   
 \end{equation}
 with
 $$
 \Sigma_{1}= \begin{pmatrix}
                                       \vspace{0.3cm}
                                       \int_0^1\sigma_{1,t}^4 \mathrm dt &  \int_0^1\sigma_{1,t}^2\sigma_{2,t}^2\rho_t^2 \mathrm dt & \int_0^1\sigma_{1,t}^3\sigma_{2,t}\rho_t \mathrm dt\\
                                        \vspace{0.3cm}
                                        \int_0^1\sigma_{1,t}^2\sigma_{2,t}^2\rho_t^2 \mathrm dt &  \int_0^1\sigma_{2,t}^4 \mathrm dt & \int_0^1\sigma_{1,t}\sigma_{2,t}^3\rho_t \mathrm dt\\
                                        \int_0^1\sigma_{1,t}^3\sigma_{2,t}\rho_t \mathrm dt & \int_0^1\sigma_{1,t}\sigma_{2,t}^3\rho_t \mathrm dt &  \int_0^1\dfrac{1}{2}\sigma_{1,t}^2\sigma_{2,t}^2(1+\rho_t^2) \mathrm dt
                                      \end{pmatrix}.
 $$
\end{thm}

 \begin{rem}
 If for some $p>2$, $\sigma_{1,t}^2$, $\sigma_{2,t}^2$ and $\sigma_{1,t}\sigma_{2,t}(1-\rho_t^2)$ $\in$ $L^p([0,1])$ and $b_n=O(n^{\frac{1}{2}-\frac{1}{p}})$, 
  the condition \eqref{Max_condition} in {\bf (MDP)} is verified.
 \end{rem}

\vspace{10pt}
 
 Let ${\mathcal H}$ be the banach space of $\mathbb{R}^3$-valued  right-continuous-left-limit non decreasing functions $\gamma$ on $[0,1]$ with $\gamma(0)=0$, equipped with the uniform norm
and the $\sigma-$field $\mathcal{B}^s$ generated by the coordinate $\{\gamma(t),0 \leqslant t \leqslant 1\}$. 

\begin{thm}
 \label{MDP2}  
 Under the conditions {\bf (MDP)} and {\bf (B)}, the sequence
 $$\dfrac{\sqrt{n}}{b_n}\left(V_{.}^n(X)-[V]_{.}\right)$$ 
 satisfies the LDP on $\mathcal H$ with speed $b_n^2$ and with the rate function given by
 \begin{equation}
 J_{mdp}(\phi)=\label{ARR}
  \left\{
  \begin{array}{lcl}
  \vspace{0.3cm}
 \displaystyle \int_0^1 \frac{1}{2}\left\langle \dot \phi(t),\Sigma_t^{-1}\cdot \dot \phi(t) \right\rangle dt\qquad if \quad \phi \in \mathcal{AC}_0([0,1])\\		
  \vspace{0.3cm}
  +\infty, \qquad otherwise,
  \end{array}
  \right.
  \end{equation}
 where 
 $$
 \Sigma_t= \begin{pmatrix}
                                       \vspace{0.3cm}
                                       \sigma_{1,t}^4 & \sigma_{1,t}^2\sigma_{2,t}^2\rho_t^2 & \sigma_{1,t}^3\sigma_{2,t}\rho_t \\
                                        \vspace{0.3cm}
                                        \sigma_{1,t}^2\sigma_{2,t}^2\rho_t^2 & \sigma_{2,t}^4 & \sigma_{1,t}\sigma_{2,t}^3\rho_t \\
                                        \sigma_{1,t}^3\sigma_{2,t}\rho_t & \sigma_{1,t}\sigma_{2,t}^3\rho_t & \dfrac{1}{2}\sigma_{1,t}^2\sigma_{2,t}^2(1+\rho_t^2) 
                                      \end{pmatrix}
 $$
 is invertible and  $\Sigma_t^{-1}$ his inverse such that
 $$
 \Sigma_t^{-1}=\frac{1}{{\rm det}(\Sigma_t)} \begin{pmatrix}
												  \vspace{0.3cm}
												  \dfrac{1}{2}\sigma_{1,t}^2\sigma_{2,t}^6(1-\rho_t^2) & \dfrac{1}{2}\sigma_{1,t}^4\sigma_{2,t}^4\rho_t^2(1-\rho_t^2) & -\sigma_{1,t}^3\sigma_{2,t}^5\rho_t(1-\rho_t^2)\\
												    \vspace{0.3cm}
												    \dfrac{1}{2}\sigma_{1,t}^4\sigma_{2,t}^4\rho_t^2(1-\rho_t^2) & \dfrac{1}{2}\sigma_{1,t}^6\sigma_{2,t}^2(1-\rho_t^2) & -\sigma_{1,t}^5\sigma_{2,t}^3\rho_t(1-\rho_t^2)\\
												    -\sigma_{1,t}^3\sigma_{2,t}^5\rho_t(1-\rho_t^2) & -\sigma_{1,t}^5\sigma_{2,t}^3\rho_t(1-\rho_t^2) &  \sigma_{1,t}^4\sigma_{2,t}^4(1-\rho_t^4)
												  \end{pmatrix},
$$												  

$$with\qquad {\rm det}(\Sigma_t)=\dfrac{1}{2}\sigma_{1,t}^6\sigma_{2,t}^6(1-\rho_t^2)^3,$$
and $\mathcal{AC}_0=\left\{ \phi:[0,1] \rightarrow \mathbb{R}^3\,\, {\it  is\, absolutely\, continuous\, with}\,\phi(0)=0\right\}.$
\end{thm}
\vspace{10pt}

Let us note once again that it is the first time the MDP is considered for the vector of (co)-volatility.\\\vspace{10pt}

In the previous results, we have imposed the boundedness of $b(t,\omega)$ which allows us to reduce quite easily the LDP and MDP of $V ^n(X)$to those of $V^n(X-Y)$ (no drift case). It is very natural to ask whether they continue to hold under a Lipchitzian condition or more generally linear growth condition of the drift $b(t,x)$, rather than the boundedness. This is the object of the following

\begin{thm}
 \label{Lipschitz}  
Let $X_t=(X_{1,t},X_{2,t})$ be given by \eqref{AR22}, with $(X_{1,0}, X_{2,0})$ bounded. We assume that the drift $b_{\ell}$ satisfies the following uniform linear growth condition: $\forall s,t \in [0,1], x,y \in \mathbb{R}^2$
\begin{equation} 
\label{yac_drift}
|b_{\ell}(t,x)-b_{\ell}(s,y)| \leqslant C[1+\|x-y\|+\eta(|t-s|)(\|x\|+\|y\|)],
\end{equation}
where $C>0$ is a constant and $\eta:[0,\infty) \to [0,\infty)$ is a continuous non-decreasing function with $\eta(0)=0$.
\begin{itemize}
 \item[(1)] Under the condition {\bf (LDP)}, the sequence $\tilde V_{\cdot}^n(X)$ satisfies the LDPs in Theorem \ref{LDP} and Theorem \ref{LDP2}. 
 \item[(2)] Under the condition {\bf (MDP)}, the sequence  $\dfrac{\sqrt{n}}{b_n}(\tilde V_{\cdot}^n(X)-[V]_{\cdot})$ satisfies the MDPs in Theorem \ref {MDP1} and Theorem \ref{MDP2}.
 \end{itemize}

\vspace{10pt}
\end{thm}
 As it can be remarked, the LDP and the MDP are established here for $\tilde V^n$ instead of $V^n$. If we conjecture that the MDP may still be valid in this case with $V^n$, we do not believe it should be the case for the LDP, and it is thus a challenging and interesting question to establish the LDP in this case for $V^n$. However for the statistical purpose, if the drift $b$ is known, the previous result is perfectly satisfactory.

 
 \section{Applications: Large deviations for the realized correlation and the realized regression coefficients }
 
 In this section we apply our results to obtain the LDP and MDP for the standard dependence measures between the two assets returns such as the realized regression coefficients up to time $1$, $\beta_{\ell, 1}=\frac{C_1}{Q_{\ell,1}}$ for $\ell =1,2$ and the realized correlation $\varrho_1=\frac{C_1}{\sqrt{Q_{1,1}Q_{2,1}}}$ which are estimated by $\beta_{\ell,1}^n(X)=\frac{C_1^n(X)}{Q^n_{\ell,1}(X)}$ and  $\varrho_1^n(X)=\frac{C_1^n(X)}{\sqrt{Q^n_{1,1}(X)Q^n_{2,1}(X)}}$ respectively. To simplify the argument, we focus in the case where $\sigma_{\ell}$ for $\ell=1,2$ are constants and we denote $\varrho:=\int_0^1\rho_tdt$. The consistency and the central limit theorem  for these estimators were already studied see for example Mancini and Gobbi \cite{Mancini2}. Up to our knowledge, however no results are known for the large and moderate deviation principle.

\subsection{Correlation coefficient}
 
 \begin{prop}
 \label{LDP3} Let for $\ell=1,2$,  $\sigma_{\ell}$ are constants and $\varrho:=\int_0^1\rho_tdt$.
 Under the conditions {\bf (LDP)} and {\bf (B)}, the sequence $\varrho_1^n(X)$ satisfies the LDP on $\mathbb{R}$ with speed n and with the good rate function given by
$$I_{ldp}^{\varrho}(u)=\inf_{\{(x,y,z)\in\mathbb{R}^3: u=\frac{z}{\sqrt{xy}}\}}I_{ldp}(x,y,z)$$
where $I_{ldp}$ is given in (\ref{LDPtaux1}).
\end{prop}
\vspace{10pt}
Once again, let us specify the rate function in the case of constant correlation.
\begin{cor}\label{LDPrhoconstant}
We suppose that for $\ell=1,2$,  $\sigma_{\ell}$ and $\rho$ are constant. Under the condition {\bf (B)}, we obtain that $\varrho_1^n(X)$ satisfies the LDP on $\mathbb{R}$ with speed n and with the good rate function given by
\begin{equation}\label{LDPrhoconstanttaux1}
{I}_{ldp}^{\rho}(u)=\left\{
  \begin{array}{lcl}
  \vspace{0.3cm}
  \log\left(\dfrac{\sqrt{1-\rho u}}{\sqrt{1-\rho^2}\sqrt{1-u^2}}\right)-1+\dfrac{\sigma_1^4 +\sigma_2^4-2\rho\sigma_1^2\sigma_2^2 u}{2\sigma_1^2\sigma_2^2(1-\rho u)}, \quad  -1<u<1\\						 \vspace{0.3cm}
  +\infty, \qquad otherwise.															  
  \end{array}
  \right.
\end{equation}
\end{cor}

As the reader can imagine from the rate function expression, it is quite a simple application of the contraction principle starting from the LDP of the realized (co)-volatility. As will be seen from the proof, in this case, the MDP is harder to establish and requires a more subtle technology: large deviations for the delta-method.

\begin{prop}
 \label{MDP3}
Let for $\ell=1,2$,  $\sigma_{\ell}$ are constants and $\varrho:=\int_0^1\rho_tdt$. Under the conditions {\bf (MDP)} and {\bf (B)}, the sequence $\dfrac{\sqrt{n}}{b_n}\left(\varrho_1^n(X)-\varrho\right)$ satisfies the LDP on $\mathbb{R}$ with speed $b_n^2$ and with the rate function given by
$$I_{mdp}^{\varrho}(u)=\inf_{\{(x,y,z)\in\mathbb{R}^3: u=\frac{z}{\sigma_1\sigma_2} - \varrho\frac{\sigma_1^2y+x\sigma_2^2}{2\sigma_1^2\sigma_2^2}\}}I_{mdp}(x,y,z)$$
where $I_{mdp}$ is given in (\ref{MDPtaux1}).
\end{prop}

\begin{cor}
We suppose that for $\ell=1,2$,  $\sigma_{\ell}$ and $\rho$ are constant. Under the condition {\bf (B)}, we obtain that $\dfrac{\sqrt{n}}{b_n}(\varrho_1^n(X)-\rho)$ satisfies the LDP on $\mathbb{R}$ with speed n and with the good rate function given for all $u\in \mathbb{R}$ by
\begin{equation}\label{MDPrhoconstanttaux1}
I_{mdp}^{\rho}(u) = \frac{2u^2}{(1-\rho^2)^2}.
\end{equation}
\end{cor}

\subsection{Regression coefficient}
The strategy initiated for the correlation coefficient is even simpler in the case of regression coefficient.

\begin{prop}
 \label{LDP4} Let for $\ell=1,2$,  $\sigma_{\ell}$ are constants .
 Under the conditions {\bf (LDP)} and {\bf (B)}, for $\ell=1$ or $2$, the sequence $\beta_{l,1}^n(X)$ satisfies the LDP on $\mathbb{R}$ with speed n and with the good rate function given by
$$I_{ldp}^{\beta_{\ell,1}}(u)=\inf_{\{(x_1,x_2,x_3)\in\mathbb{R}^3: u=\frac{x_3}{x_{\ell}}\}}I_{ldp}(x_1,x_2,x_3)$$
where $I_{ldp}$ is given in (\ref{LDPtaux1}).
\end{prop}

Once again, this Proposition is a simple application of the contraction principle. Let us specify the rate function in the case of constant correlation.

\begin{cor}We suppose that for $\ell=1,2$,  $\sigma_{\ell}$ and $\rho$ are constant. Under the condition {\bf (B)}, we obtain that $\beta_{l,1}^n(X)$ satisfies the LDP on $\mathbb{R}$ with speed n and with the good rate function given for $\iota=1,2$ with $\ell\not=\iota$ and for all $u\in \mathbb{R}$ by
\begin{equation}\label{LDPrhoconstanttaux12}
I_{ldp}^{\beta_l}(u)=\frac{1}{2} \log\left(1+\frac{(\sigma_{\ell} u-\rho \sigma_{\iota})^2}{\sigma_{\iota}^2(1-\rho^2)}\right).
\end{equation}
\end{cor}

We may also consider the MDP.

\begin{prop}
 \label{MDP4} Let for $\ell=1,2$,  $\sigma_{\ell}$ are constants and $\varrho:=\int_0^1\rho_tdt$. Under the conditions {\bf (MDP)} and {\bf (B)} and for $\ell,\iota\in\{1,2\}$ with $\ell\not=\iota$, the sequence $\frac{\sqrt{n}}{b_n}(\beta_{\ell,1}^n(X)-\varrho\frac{\sigma_{\iota}}{\sigma_{\ell}})$ satisfies the LDP on $\mathbb{R}$ with speed $b_n^2$ and with the rate function given by
$$I_{mdp}^{\beta_{\ell,1}}(u)=\inf_{\{(x,y,z)\in\mathbb{R}^3: u=\frac{z}{\sigma_{\ell}^2}-\varrho\frac{\sigma_{\iota}}{\sigma_{\ell}^3}x\}}I_{mdp}(x,y,z)$$
where $I_{mdp}$ is given in (\ref{MDPtaux1}).
\end{prop}

\begin{cor}
We suppose that for $\ell=1,2$,  $\sigma_{\ell}$ and $\rho$ are constant. Under the condition {\bf (B)} and for $\ell,\iota\in\{1,2\}$ with $\ell\not=\iota$, we obtain that  $\frac{\sqrt{n}}{b_n}(\beta_{\ell,1}^n(X)-\rho\frac{\sigma_{\iota}}{\sigma_{\ell}})$ satisfies the LDP on $\mathbb{R}$ with speed n and with the good rate function given for all $u\in \mathbb{R}$ by
\begin{equation}
\label{MDPrhoconstanttaux2}
I_{mdp}^{\beta^c_{\ell,1}}(u) =\frac{2\sigma_{\ell}^2 u^2}{\sigma_{\iota}^2(1-\rho^2)}.
\end{equation}

\end{cor}

\section{Proof}

Let us say a few words on our strategy of proof. As the reader may have guessed, one of the important step is first to consider the no-drift case, where we have to deal with non homogenous quadratic forms of Gaussian processes (in the vector case). In these non essentially smooth case (in the terminology of G\"artner-Ellis), we will use (after some technical approximations) powerful recent results of Najim \cite{Najim1}. In a second step, we see how to reduce the general case to the no-drift case.

\subsection{Proof of  Theorem \ref{LDP}}
\begin{lem}\label{lem_1}
 If $(\xi,\xi^{'})$ are independent centered Gaussian random vector with covariance
 $$
\begin{pmatrix}
 1 & c\\
 c & 1
\end{pmatrix}
,-1<c<1.
$$
Then for all $(\lambda_1,\lambda_2,\lambda_3) \in \mathbb{R}^3$ 
$$\displaystyle\log\mathbb{E}\exp\left(\lambda_1\xi^2+\lambda_2\xi^{'2}+\lambda_3\xi\xi^{'}\right)=P_c(\lambda_1,\lambda_2,\lambda_3),$$
where the function $P_c$ is given in in (\ref{ASSA}).
\end{lem}

\subsection*{Proof :} Elementary.

\begin{lem}\label{lem_2}
Let $X_t=(X_{1,t},X_{2,t})$ given (\ref{AR}) and $Y_t=(Y_{1,t},Y_{2,t})$ where for $\ell =1,2\quad$$Y_{\ell,t}:=\int_0^tb_{\ell}(t,\omega)dt$. We have for every $\lambda \in \mathbb{R}^3$
\begin{equation*}
\Lambda_n(\lambda):=\dfrac{1}{n}\log\mathbb{E}\left(\exp\left(n\left\langle \lambda,V_1^n (X-Y)\right\rangle\right)\right)\leqslant\Lambda(\lambda):=\int_0^1 P_{\rho_t}(\lambda_1\sigma_{1,t}^2,\lambda_2\sigma_{2,t}^2,\lambda_3\sigma_{1,t}\sigma_{2,t}) \mathrm dt,
\end{equation*}
where the function $P_c$ is given in in (\ref{ASSA}), and
 $$
 \lim_{n \to \infty} \Lambda_n(\lambda) =\Lambda(\lambda).
 $$ 
 \end{lem}                        
                               
\subsection*{Proof :}
For $\ell =1,2$, we have
$$Q_{\ell,t}^n(X-Y)=\sum_{k=1}^{[nt]}a_{\ell,k}\xi_{\ell,k}^2\qquad{\rm and}\quad C_{t}^n(X-Y)=\sum_{k=1}^{[nt]}\sqrt{a_{1,k}} \sqrt{a_{2,k}}\xi_{1,k}\xi_{2,k}
$$
where
\begin{equation}\label{defxi}
\xi_{\ell,k}:=\dfrac{\int_{t_{k-1}}^{t_k}\sigma_{\ell,s}\mathrm dB_{\ell,s}}{\sqrt{a_{\ell,k}}}\quad {\rm and} \quad a_{\ell,k}:=\int_{t_{k-1}^n}^{t_k^n}\sigma_{\ell,s}^2\mathrm ds.
\end{equation}

Obviously $((\xi_{1,k},\xi_{2,k}))_{k=1 \cdots n}$ are independent centered Gaussian random vector with covariance matrix
$$
\begin{pmatrix}
\vspace{0.3cm}
 1 & c_k^n \\
c_k^n & 1
\end{pmatrix}
$$ 
where  
\begin{equation}\label{defc}
c_k^n:=\dfrac{ \vartheta_k^n}{\sqrt{a_{1,k}}\sqrt{a_{2,k}}}\qquad{\rm and} \quad\vartheta_k^n:=\int_{t_{k-1}^n}^{t_k^n}\sigma_{1,s}\sigma_{2,s}\rho_s \mathrm ds.
\end{equation} 

We use the lemma \ref{lem_1} and the martingale convergence theorem (or the classical Lebesgue derivation theorem) to get the final assertions (see for example \cite[p.204]{Djellout1} for details).
\vspace{10pt}

\subsection*{Proof Theorem \ref{LDP}}$\quad$

\vspace{10pt}
{\rm (1)} It is contained in lemma \ref{lem_2}.

\vspace{10pt}

{\rm (2)} We shall prove it in three steps.
\vspace{10pt}

\vspace{10pt}

{\underline{\it Part 1}.} At first, we consider that the drift $b_{\ell}=0$. In this case $V_t^n(X)=V_t^n(X-Y)$. We recall that since $B_{2,t}=\rho_tdB_{1,t} + \sqrt{1-\rho_t^2}dB_{3,t}$, we may rewrite (\ref{AR}) as  
 \begin{equation}
\label{AR2}
\vspace{1ex}
\left\{
\begin{array}[c]{ccccc}
dX_{1,t} & = & \sigma_{1,t}dB_{1,t} \\
dX_{2,t} & = & \sigma_{2,t}(\rho_tdB_{1,t} & + & \sqrt{1-\rho_t^2}dB_{3,t})
\end{array}
\right.
\end{equation}

Using the approximation Lemma in \cite{Demzei}, we shall prove that
$$V_1^n(X-Y)=\left( \sum_{k=1}^{n}\left(\Delta_k^nX_1\right)^2, \sum_{k=1}^{n}\left(\Delta_k^nX_2\right)^2, \sum_{k=1}^{n}\left(\Delta_k^nX_1\right)\left(\Delta_k^nX_2\right)       \right)^T$$
will satisfy the same LDPs as
$$W_{1}^n:=\begin{pmatrix}\displaystyle
                                       \vspace{0.3cm}
                                       \dfrac{1}{n}\sum_{k=1}^{n}\sigma_{1,\frac{k-1}{n}}^2N_{1,k}^2\\
                                        \vspace{0.3cm}\displaystyle
                                       \dfrac{1}{n}\sum_{k=1}^{n} \left(\sigma_{2,\frac{k-1}{n}}\rho_{\frac{k-1}{n}}N_{1,k} + \sigma_{2,\frac{k-1}{n}}\sqrt{1-\rho_{\frac{k-1}{n}}^2} N_{3,k} \right)^2\\ \displaystyle
                                       \dfrac{1}{n}\sum_{k=1}^{n} \sigma_{1,\frac{k-1}{n}}N_{1,k}\left(\sigma_{2,\frac{k-1}{n}}\rho_{\frac{k-1}{n}}N_{1,k} + \sigma_{2,\frac{k-1}{n}}\sqrt{1-\rho_{\frac{k-1}{n}}^2} N_{3,k} \right)\\
                                       \end{pmatrix},  
$$
 where $N_{\ell,k}:=\int_{t_{k-1}^n}^{t_k^n}\sqrt{n} \mathrm dB_{\ell,s}$, for $\ell=1,3.$

Let us first focus on the LDP of $W_1^n$. We will use Najim result (see Lemma \ref{Najim_Jamal_1}) to prove that. 

It is easy to see that $W_1^n$ can be reritten as
$$W_1^n=\dfrac{1}{n}\sum_{k=1}^{n}F\left(\frac{k-1}{n}\right)Z_k $$
where
\begin{equation}\label{defF}
F\left(\frac{k}{n}\right)= \begin{pmatrix}
                                       \vspace{0.3cm}
                                       f_1(\frac{k}{n})\\
                                        \vspace{0.3cm}
                                       f_2(\frac{k}{n})\\
                                       f_3(\frac{k}{n})\\
                                       \end{pmatrix}
    = \begin{pmatrix}
                                       \vspace{0.3cm}
                                       \sigma_{1,\frac{k}{n}}^2 & 0 & 0 \\
                                        \vspace{0.3cm}
                                       \sigma_{2,\frac{k}{n}}^2\rho_{\frac{k}{n}}^2 & \sigma_{2,\frac{k}{n}}^2(1-\rho_{\frac{k}{n}}^2) & 2\sigma_{2,\frac{k}{n}}^2\rho_{\frac{k}{n}}\sqrt{1-\rho_{\frac{k}{n}}^2}\\
                                       \sigma_{1,\frac{k}{n}}\sigma_{2,\frac{k}{n}}\rho_{\frac{k}{n}} & 0 & \sigma_{1,\frac{k}{n}}\sigma_{2,\frac{k}{n}}\sqrt{1-\rho_{\frac{k}{n}}^2}
                                       \end{pmatrix}
\end{equation}
and
\begin{equation}\label{Z}
Z_j = \left(N_{1,j}^2,\,\, N_{3,j}^2,\,\, N_{1,j} N_{3,j}\right)^T.
\end{equation}

Obviously $(N_{1,k},N_{3,k})_{k=1 \cdots n}$ are independent centered Gaussian random vector with identity covariance matrix.

For the LDP of $W_1^n$ we will use Lemma 5.1, in the case where ${\mathcal X}:=[0,1]$ and $R(dx)$ is the Lebesgue measure on $[0,1]$ and $x_i^n:=i/n$. One can check that, in this situation, Assumptions (N-2) hold true. The random variables $(Z_k)_{k=1,\cdots,n}$ are independent and identically distributed. By the definition of $Z_k$, the Assumptions (N-1) hold true also.  

So $W_{1}^n$ satisfies the LDP on $\mathbb{R}^3$ with speed $n$ and with the good rate function given by  for all $x\in \mathbb{R}^3$
$$I(x)=\sup_{\lambda\in\mathbb{R}^3}\left(\left\langle \lambda,x \right\rangle - \int_0^1L\left(\sum_{j=1}^{3}\lambda_i\cdot f_i(t)\right)dt\right),$$
with
$$
\sum_{i=1}^{3}\lambda_i f_i(t) = \begin{pmatrix}
                                       \vspace{0.3cm}
                                       \lambda_1\sigma_{1,t}^2+\lambda_2\sigma_{2,t}^2\rho_{t}^2+\lambda_3\sigma_{1,t}\sigma_{2,t}\rho_{t}\\
                                        \vspace{0.3cm}
                                      \lambda_2\sigma_{2,t}^2(1-\rho_{t}^2)
\\
                                       2\lambda_2\sigma_{2,t}^2\rho_{t}\sqrt{1-\rho_{t}^2}+\lambda_3\sigma_{1,t}\sigma_{2,t}\sqrt{1-\rho_{t}^2}
                                       \end{pmatrix}^T,
$$
and for $\lambda\in \mathbb{R}^3$ 
$$L(\lambda):=\log\mathbb{E}\exp{\langle\lambda, Z_1\rangle}=P_0(\lambda_1,\lambda_2,\lambda_3),$$
where $P_0$ is given in (\ref{ASSA}). In this cas it takes a simpler form wich we recall here:
$$P_0(\lambda_1,\lambda_2,\lambda_3)=-\frac{1}{2}\log\left((1-2\lambda_1)(1-2\lambda_2)-\lambda_3^2\right).$$

An easy calculation gives us that
$$\int_0^1P_0\left(\sum_{j=1}^{3}\lambda_i\cdot f_i(t)\right)dt=\int_0^1P_{\rho_t}\left(\lambda_1\sigma^2_{1,t}, \lambda_2\sigma^2_{2,t},\lambda_3\sigma_{1,t}\sigma_{2,t}\right)dt,$$
so 
$$I(x)=I_{ldp}(x),$$
where $I_{ldp}$ is given in (\ref{LDPtaux1}).

\vspace{10pt}

{\underline{\it Part 2}.} Now we shall prove that $V_1^n$  and $W_1^n$ satisfy the same LDPs, by means of the approximation Lemma in \cite{Djellout1}.  We have to prove that 
$$ V_1^n(X-Y) - W_1^n\superexpldp 0.
$$

We do this element by element. We will only consider one element, the other terms can be dealt with in the same way. We have to prove that  for $q=1,2,3$
\begin{equation}\label{negR}
{\mathcal R}_{q,1}^n\superexpldp 0,
\end{equation}
where
\begin{equation}\label{calR1}
{\mathcal R}_{1,t}^n:=\sum_{k=1}^{[nt]}\left(\Delta_k^nX_1\right)^2-\frac{1}{n}\sum_{k=1}^{[nt]}\sigma_{1,\frac{k-1}{n}}^2N_{1,k}^2,
\end{equation}
\begin{equation}\label{calR2}
{\mathcal R}_{2,t}^n:=\sum_{k=1}^{[nt]}(\Delta_k^nX_2)^2 - \frac{1}{n}\sum_{k=1}^{[nt]}\biggl(\sigma_{2,\frac{k-1}{n}}\rho_{\frac{k-1}{n}}N_{1,k}+\sigma_{2,\frac{k-1}{n}}\sqrt{1-\rho_{\frac{k-1}{n}}^2} N_{3,k}\biggr)^2,
\end{equation}
and
\begin{equation}\label{calR3}
{\mathcal R}_{3,t}^n:=\sum_{k=1}^{[nt]}\Delta_k^nX_1 \Delta_k^nX_2- \frac{1}{n}\sum_{k=1}^{[nt]}\sigma_{1,\frac{k-1}{n}}N_{1,k}\biggl(\sigma_{2,\frac{k-1}{n}}\rho_{\frac{k-1}{n}}N_{1,k}+\sigma_{2,\frac{k-1}{n}}\sqrt{1-\rho_{\frac{k-1}{n}}^2} N_{3,k}\biggl).
\end{equation}

At first, we start the negligibility (\ref{negR}) with the quantity ${\mathcal R}_{1,1}^n$ which can be rewritten as
 $$
{\mathcal R}_{1,1}^n=\sum_{k=1}^{n}R_{-,k} R_{+,k},$$
with $R_{\pm,k}:=\int_{t_{k-1}^n}^{t_k^n}(\sigma_{1,s}\pm\sigma_{1,\frac{k-1}{n}}) \mathrm dB_{1,s},$ where $((R_{-,k},R_{+,k}))_{k=1 \cdots n}$ are independent centered Gaussian random vector with covariance
$$
 \begin{pmatrix}
 \vspace{0.3cm}
  \varepsilon_{-,k}^n& \eta_k^n  \\
   \eta_k^n &  \varepsilon_{+,k}^n
 \end{pmatrix}				 
$$
where 
\begin{equation}\label{defeta}
\varepsilon_{\pm,k}^n=\int_{t_{k-1}^n}^{t_k^n}\left(\sigma_{1,s}\pm\sigma_{1,\frac{k-1}{n}}\right)^2 \mathrm ds\qquad{\rm and}\qquad  \eta_k^n=\int_{t_{k-1}^n}^{t_k^n}\left(\sigma_{1,s}^2-\sigma_{1,\frac{k-1}{n}}^2\right) \mathrm ds
\end{equation}

So by Chebyshev's inequality, we have for all $r,\lambda > 0,$
\begin{equation} \label{LDp11}
\frac{1}{n}\log\mathbb{P}\left({\mathcal R}_{1,1}^n> r\right)
\leqslant -r\lambda + \frac{1}{n}\log\mathbb{E}\exp\left(n\lambda {\mathcal R}_{1,1}^n\right)\\
\end{equation}

A simple calculation gives us
\begin{eqnarray}\label{K}
\frac{1}{n}\log\mathbb{E}\exp\left(n\lambda {\mathcal R}_{1,1}^n\right)&=&\frac{1}{n}\sum_{k=1}^{n}\log\mathbb{E}\exp\left(n\lambda R_{+,k} R_{-,k}\right)\nonumber\\
 &=& -\frac{1}{2n}\sum_{k=1}^{n}\log\left[\dfrac{\varepsilon_{+,k}^n\varepsilon_{-,k}^n-\left(n\lambda(\varepsilon_{+,k}^n\varepsilon_{-,k}^n-(\eta_k^n)^2)+\eta_k^n\right)^2}{\varepsilon_{+,k}^n\varepsilon_{-,k}^n-(\eta_k^n)^2} \right]\nonumber\\
&=& -\frac{1}{2n}\sum_{k=1}^{n}\log\left[ 1-n^2\lambda^2(\varepsilon_{+,k}^n\varepsilon_{-,k}^n-(\eta_k^n)^2)-n\lambda \eta_k^n\right]\nonumber\\
&=&\int_0^1 K(f_n(t)) \mathrm dt
\end{eqnarray}
where $K$ is given by
$$
K(\lambda):=\left\{
  \begin{array}{lcl}
  \vspace{0.3cm}
  -\dfrac{1}{2} \log\left(1-2\lambda\right)  \quad if \quad \lambda<\frac{1}{2}\\
   \vspace{0.3cm}
  +\infty, \qquad otherwise,
  \end{array}
  \right.
$$
and
$$f_n(t)=\sum_{k=1}^{n}1_{(t_{k}^n-t_{k-1}^n)}(t)\left[\lambda^2\left(\frac{\varepsilon_{+,k}^n}{t_{k}^n-t_{k-1}^n}\frac{\varepsilon_{-,k}^n}{t_{k}^n-t_{k-1}^n}-\left(\frac{\eta_k^n}{t_{k}^n-t_{k-1}^n}\right)^2\right)+2\lambda\left(\frac{\eta_k^n}{t_{k}^n-t_{k-1}^n}\right)\right]$$
where $\varepsilon_{\pm,k}^n$ and $\eta_k^n$ are given in (\ref{defeta}).

By the continuity condition of the assumption {\bf (LDP)} and the classical Lebesgue derivation theorem, we have that
$$f_n(t) \longrightarrow 0 \qquad{\rm as}\quad n\rightarrow \infty.$$

By the classical Lebesgue derivation theorem we have that the right hand of the equality (\ref{K}) goes to $0$.

Letting $n$ goes to infinity and than $\lambda$ goes to infinity in (\ref{LDp11}), we obtain that
$$
\lim_{n\rightarrow \infty}\frac{1}{n}\log\mathbb{P}\left({\mathcal R}_{1,1}^n> r\right)=-\infty.
$$

Doing the same things with $-{\mathcal R}_{1,1}^n$, we obtain (\ref{negR}) for ${\mathcal R}_{1,1}^n$.
\vspace{10pt}

Now we shall prove \eqref{negR} with ${\mathcal R}_{2,1}^n$. We have 
$${\mathcal R}_{2,1}^n=\sum_{k=1}^{n}E_{-,k} E_{+,k} + \sum_{k=1}^{n}E_{-,k} B_{+,k} + \sum_{k=1}^{n}E_{+,k} B_{-,k} + \sum_{k=1}^{n}B_{-,k} B_{+,k},$$
where
$$E_{\pm,k}:=\int_{t_{k-1}^n}^{t_k^n}(\sigma_{2,s}\rho_{s}\pm\sigma_{2,\frac{k-1}{n}}\rho_{\frac{k-1}{n}}) \mathrm dB_{1,s},$$
and
$$B_{\pm,k}:=\int_{t_{k-1}^n}^{t_k^n}(\sigma_{2,s}\sqrt{1-\rho_{s}^2}\pm\sigma_{2,\frac{k-1}{n}}\sqrt{1-\rho_{\frac{k-1}{n}}^2}) \mathrm dB_{3,s}.$$

where $(E_{-,k}, E_{+k})$,$(E_{-,k},B_{+,k})$,$(E_{+,k},B_{-,k})$,$(B_{-,k},B_{+,k})$, ${k=1 \cdots n}$ are four independent centered Gaussian random vectors with covariances respectively given by

{\begin{center}
$$
 \begin{pmatrix}
 \vspace{0.3cm}
 \int_{t_{k-1}^n}^{t_k^n}(\sigma_{2,s}\rho_{s}-\sigma_{2,\frac{k-1}{n}}\rho_{\frac{k-1}{n}})^2 \mathrm ds & \int_{t_{k-1}^n}^{t_k^n}(\sigma_{2,s}^2\rho_{s}^2-\sigma_{2,\frac{k-1}{n}}^2\rho_{\frac{k-1}{n}}^2) \mathrm ds \\
 \int_{t_{k-1}^n}^{t_k^n}(\sigma_{2,s}^2\rho_{s}^2-\sigma_{2,\frac{k-1}{n}}^2\rho_{\frac{k-1}{n}}^2) \mathrm ds &  \int_{t_{k-1}^n}^{t_k^n}(\sigma_{2,s}\rho_{s}+\sigma_{2,\frac{k-1}{n}}\rho_{\frac{k-1}{n}})^2 \mathrm ds
 \end{pmatrix}
 $$
 \end{center}
 
 and
 
{\begin{center}
$$
 \begin{pmatrix}
 \vspace{0.3cm}
 \int_{t_{k-1}^n}^{t_k^n}(\sigma_{2,s}\rho_{s}-\sigma_{2,\frac{k-1}{n}}\rho_{\frac{k-1}{n}})^2 \mathrm ds & 0 \\
 0 &  \int_{t_{k-1}^n}^{t_k^n}(\sigma_{2,s}\sqrt{1-\rho_{s}^2}+\sigma_{2,\frac{k-1}{n}}\sqrt{1-\rho_{\frac{k-1}{n}}^2})^2 \mathrm ds
 \end{pmatrix}
 $$
  \end{center}				 

and

{\begin{center}
$$
 \begin{pmatrix}
 \vspace{0.3cm}
 \int_{t_{k-1}^n}^{t_k^n}(\sigma_{2,s}\rho_{s}+\sigma_{2,\frac{k-1}{n}}\rho_{\frac{k-1}{n}})^2 \mathrm ds & 0 \\
 0 &  \int_{t_{k-1}^n}^{t_k^n}(\sigma_{2,s}\sqrt{1-\rho_{s}^2}-\sigma_{2,\frac{k-1}{n}}\sqrt{1-\rho_{\frac{k-1}{n}}^2})^2 \mathrm ds
 \end{pmatrix}
 $$
 \end{center}			 

and

{\begin{center}
$$
 \begin{pmatrix}
 \vspace{0.3cm}
  \int_{t_{k-1}^n}^{t_k^n}(\sigma_{2,s}\sqrt{1-\rho_{s}^2}-\sigma_{2,\frac{k-1}{n}}\sqrt{1-\rho_{\frac{k-1}{n}}^2})^2 \mathrm ds & \int_{t_{k-1}^n}^{t_k^n}(\sigma_{2,s}^2(1-\rho_{s}^2)-\sigma_{2,\frac{k-1}{n}}^2(1-\rho_{\frac{k-1}{n}}^2)) \mathrm ds \\
  \int_{t_{k-1}^n}^{t_k^n}(\sigma_{2,s}^2(1-\rho_{s}^2)-\sigma_{2,\frac{k-1}{n}}^2(1-\rho_{\frac{k-1}{n}}^2)) \mathrm ds &  \int_{t_{k-1}^n}^{t_k^n}(\sigma_{2,s}\sqrt{1-\rho_{s}^2}+\sigma_{2,\frac{k-1}{n}}\sqrt{1-\rho_{\frac{k-1}{n}}^2})^2 \mathrm ds.
 \end{pmatrix}
 $$	
  \end{center}
  
So \eqref{negR} for ${\mathcal R}_{2,1}^n$ is deduced if
  
  $$\sum_{k=1}^{n}E_{-k} E_{+,k}\superexpldp 0,\qquad \sum_{k=1}^{n}E_{-,k} B_{+,k} \superexpldp 0,$$
  $$\sum_{k=1}^{n}E_{+,k} B_{-,k} \superexpldp 0,\qquad \sum_{k=1}^{n}B_{-,k} B_{+,k} \superexpldp 0.$$
  
Each convergence  is deduced by the same calculations as for  $(R_{-,k},R_{+,k})_{k=1 \cdots n}$.
\vspace{10pt}

Now we shall prove \eqref{negR} with ${\mathcal R}_{3,1}^n$. We have 
\begin{eqnarray*}
{\mathcal R}_{3,1}^n&=&\sum_{k=1}^{n}R_{-,k}E_{-,k}+ \sum_{k=1}^{n} R_{-,k} B_{-,k} + \sum_{k=1}^{n} R_{-,k}\left(\int_{t_{k-1}^n}^{t_k^n} \sigma_{2,\frac{k-1}{n}}\rho_{\frac{k-1}{n}} \mathrm dB_{1,s}\right)\\
&&+ \sum_{k=1}^{n} R_{-,k} \left(\int_{t_{k-1}^n}^{t_k^n}\sigma_{2,\frac{k-1}{n}}\sqrt{1-\rho_{\frac{k-1}{n}}^2} \mathrm dB_{3,s}\right)+ \sum_{k=1}^{n} E_{-,k} \left(\int_{t_{k-1}^n}^{t_k^n}\sigma_{1,\frac{k-1}{n}} \mathrm dB_{1,s}\right)\\
&&+ \sum_{k=1}^{n} B_{-,k}\left(\int_{t_{k-1}^n}^{t_k^n}\sigma_{1,\frac{k-1}{n}} \mathrm dB_{1,s}\right),
\end{eqnarray*}
where we have used the same notation as before. By the same calculations used to prove \eqref{negR}  for ${\mathcal R}_{1,1}^n$  and ${\mathcal R}_{2,1}^n$, we obtain \eqref{negR}  for ${\mathcal R}_{3,1}^n$.
\vspace{10pt} 

Then $V_1^n(X-Y)$ and $W_1^n$ satisfy the same LDP.
\vspace{10pt}

{\underline{\it Part 3}.} We will prove that $V_1^n(X)$ and $V_1^n(X-Y)$ satisfy the same LDP. We need to prove that
$$\left\|V_1^n(X)-V_1^n(X-Y)\right\|\superexpldp 0.$$

This will be done element by element :  for $\ell =1,2$
\begin{equation}\label{sansb} Q_{\ell,1}^n(X)-Q_{\ell,1}^n(X-Y)\superexpldp 0\quad{\rm and}\quad C_1^n(X)-C_1^n(X-Y) \superexpldp 0.\end{equation}

We have 
$$\left|Q_{\ell,1}^n(X)-Q_{\ell,1}^n(X-Y)\right|\le \varepsilon(n)Q_{\ell,1}^n(X-Y)+\left(1+\frac{1}{\varepsilon(n)}\right)Z^n_{\ell},$$
and
$$\left|C_{1}^n(X)-C_{1}^n(X-Y)\right|\le \varepsilon(n)\left(Q_{1,1}^n(X-Y)+Q_{2,1}^n(X-Y)\right)+\left(\frac 12+\frac{1}{\varepsilon(n)}\right)\left(Z^n_{1}+Z^n_{2}\right).$$
with
\begin{equation*}
Z_{\ell,n}=\sum_{k=1}^{n}\left(\int_{t_{k-1}^n}^{t_k^n}b_{\ell}(t,\omega)\mathrm dt \right)^2\le \frac{\|b_{\ell}\|_{\infty}^2}{n}.
\end{equation*}
We chose $\varepsilon(n)$ such that $n\varepsilon(n)\rightarrow\infty$, so (\ref{sansb}) follows from the LDP of $Q_{\ell,1}(X)$, $C_{1}^n(X)$ and the estimations above.

\vspace{10pt}

\subsection{Proof  of Corollary \ref{LDPconstant}}$\quad$
\vspace{10pt}

From Theorem \ref{LDP}, we obtain that $V_{1}^n(X-Y)$ satisfies the LDP on $\mathbb{R}^3$ with speed $n$ and with the good rate function given by for all  $x\in \mathbb{R}^3$
$$I^{V}_{ldp}(x)=\sup_{\lambda\in\mathbb{R}^3}\left(\left\langle \lambda,x \right\rangle - P_{\rho}(\sigma_1^2\lambda_1,\sigma^2_2\lambda_2,\sigma_1\sigma_2\lambda_3)\right)=P_{\rho}^*\left(\frac{\lambda_1}{\sigma_1^2}, \frac{\lambda_2}{\sigma_2^2},\frac{\lambda_3}{\sigma_1\sigma_2}\right).
 $$
where $P_{\rho}$ and $P_{\rho}^*$ are given in (\ref{ASSA}) and (\ref{Legendre}) respectively. So we get the expression of $I^{V}_{ldp}$  given in (\ref{LDPtauxconstantV}).

The Legendre transformation of $P_c$ is defined by
$$P_{\rho}^*(x) :=\sup_{\lambda\in\mathbb{R}^3}\left(\left\langle \lambda,x \right\rangle - P_{\rho}(\lambda)\right).$$

The function $\lambda\to \left\langle \lambda,x \right\rangle - P_{\rho}(\lambda)$ reaches the supremum at the point $\lambda^*=(\lambda_1^*,\lambda_2^*,\lambda_3^*)$ such as
$$
\left\{
  \begin{array}{cc}
  \vspace{0.3cm}
  \lambda_1^*=\frac{1}{2}\dfrac{x_1x_2-(1-\rho^2)x_2-x_3^2}{2(1-\rho^2)(x_1x_2-x_3^2)},\\	
   \vspace{0.3cm}
     \lambda_2^*=\frac{1}{2}\dfrac{x_1x_2-(1-\rho^2)x_1-x_3^2}{2(1-\rho^2)(x_1x_2-x_3^2)},\\															 
     \vspace{0.3cm}
    \lambda_3^*=\dfrac{x_3^2\rho-x_1x_2\rho+(1-\rho^2)x_3}{(1-\rho^2)(x_1x_2-x_3^2)}.															  
  \end{array}
  \right.
$$

So we get the expression of the Legendre transformation $P_{\rho}^*$ given in (\ref{Legendre}).


\subsection{Proof of Theorem \ref{LDP2}}$\qquad$
\vspace{10pt}

Now we shall prove the Theorem \ref{LDP2} in two steps.  

{\underline{Step 1}. } We start by proving that the LDP holds for 
$$W_t^n=\dfrac{1}{n}\sum_{k1}^{[nt]}F\left(\frac{k-1}{n}\right)Z_k, $$
where $F$ is given in (\ref{defF}) and $Z_k$ is given in (\ref{Z}). 

This result come from an application of LDP of Lemma 5.2 derived in the case where ${\mathcal X}:=[0,1]$ and $R(dx)$ is the Lebesgue measure on $[0,1]$ and $x_i^n:=i/n$. One can check that, in this situation,
Assumptions (N-2) and (N-3) hold true.

The random variables $(Z_k)_{k=1,\cdots,n}$ are independent and identically distributed. And we will apply the Lemma with the random variables
$\displaystyle Z_k^n:=F\left(\frac{k-1}{n}\right)Z_k$. The law of $Z_k^n$ depends on the position $x_i^n:=i/n$  This type of model was partially examined by Najim see section 2.4.2 in \cite{Najim2}.

By the definition of $Z_k^n$, the Assumptions (N-1) and (N-4)  hold true.  

Finally, we just need to verify that if $x_i^n$ and $x_j^n$ are close then so are ${\mathcal L}(Z^n_i)$ and  ${\mathcal L}(Z^n_j)$ for the following Wasserstein type distance between probability measures:
$$d_{OW}(P,Q)= \inf_{\eta \in M(P,Q)} \inf \biggl\{a>0; \int_{\mathbb{R}^3\times\mathbb{R}^3}\tau\left(\frac{z-z'}{a}\right) \mathrm \eta(dzdz') \leqslant 1 \biggr\},$$
where $\eta$ is a probability with given marginals $P$ and $Q$ and $\eta(z) = e^{|z|}-1$.

In fact, consider the random variables $Y=F(x) Z$ and $\tilde Y=F(x') Z$, since $F$ is continuous
$$\mathbb E \tau\left(\frac{Y-\tilde Y}{\epsilon}\right)= \mathbb E \tau\left(\frac{(F(x)-F(x'))Z}{\epsilon}\right)\le 1$$
for $x'$ close to $x$. Thus $d_{OW}({\mathcal L}(Z^n_i) ,{\mathcal L}(Z^n_j))\le \epsilon$. This gives the Assumption (N-5).

So we deduce that the sequence $W_{\cdot}^n$ satisfies the LDP on $\mathcal BV$ with speed $n$ and with the rate function $J_{ldp}$ given by
 \begin{equation*}\label{LDPtaux2}
 J_{ldp}(f)=\int_0^1\Lambda^*_t(f_a^{'}(t)) \mathrm dt + \int_0^1\hbar_t(f_s^{'}(t)) \mathrm d\theta(t).
 \end{equation*}
 where  for all $z \in \mathbb{R}^3$ and all $t\in [0,1]$
 $$\Lambda^*_t(z)=\sup_{\lambda \in \mathbb{R}^3} \left(\left\langle \lambda,z \right\rangle - \Lambda_t(\lambda) \right),$$
 with  
 $$\Lambda_t(\lambda)=\log \int_{\mathbb{R}^3}e^{\langle \lambda,z\rangle}P(t,dz)=P_{\rho_t}(\lambda_1\sigma_{1,t}^2,\lambda_2\sigma_{2,t}^2,\lambda_3\sigma_{1,t}\sigma_{2,t}),$$
 so $\Lambda^*_t$ coincide with $P_{\rho_t}^*$ given in Theorem \ref{LDP2}. 
 
 And $\theta$ is any real-valued nonnegative measure with respect to which $\mu_s^f$ is absolutely continuous and $f_s'=d\mu_s^f/d\theta$ and for all $z \in \mathbb{R}^3$ and all $t\in [0,1]$ the recession function $\hbar_t$ of $\Lambda^*_t$ defined by $\hbar_t(z)=\sup\{\left\langle \lambda,z \right\rangle, \lambda \in \mathcal D_{\Lambda_t}\}$ with $\mathcal D_{\Lambda_t}=\{\lambda \in \mathbb{R}^3, \Lambda_t(\lambda)<\infty\}=\{\lambda \in \mathbb{R}^3, P_{\rho_t}(\lambda_1\sigma_{1,t}^2,\lambda_2\sigma_{2,t}^2,\lambda_3\sigma_{1,t}\sigma_{2,t}),<\infty\}$.

The recession function  $\alpha$ of $P^*_c$ , see Theorem 13.3 in \cite{Rock} is given by
$$\alpha(z):=\lim_{h\rightarrow \infty}\frac{P^*_c(hz)}{h}=\dfrac{z_1+z_2-2cz_3}{2(1-c^2)}1_{[z_1>0,z_2>0,z^2_3<z_1z_2]}.
$$
Using this expression, we obtain the rate function given in the Theorem \ref{LDP2}.
\vspace{10pt}

{\underline{Step 2}. Now we have to prove that
\begin{equation*}
\sup_{t\in [0,1]}\left\| V_t^n(X-Y)-W_t^n\right\|\superexpldp 0.
\end{equation*}

To do that, we have to prove that  for $q=1,2,3$
\begin{equation}\label{negR123}
\sup_{t\in [0,1]}\left|{\mathcal R}_{q,t}^n\right|\superexpldp 0,
\end{equation}
where the definition of ${\mathcal R}_{q,t}^n$ arge given in (\ref{calR1}), (\ref{calR2}) and (\ref{calR3}) for $q=1,2,3$ respectively.

We start by proving (\ref{negR123}) for $q=1$, the other terms  for $q=2,3$ follow the same line of proof.

We remark that  $({\mathcal R}_{1,t}^n-\mathbb{E}({\mathcal R}_{1,t}^n))$ is a (${\mathcal F}_{[nt]/n}$)-martingale. Then 
$$\exp(\lambda n\left[{\mathcal R}_{1,t}^n-\mathbb{E}({\mathcal R}_{1,t}^n)\right])$$
is a sub-martingale. By the maximal inequality, we have for any $r, \lambda >0$,

\begin{eqnarray*}\mathbb{P}\left(\sup_{t\in [0,1]}\left[{\mathcal R}_{1,t}^n-\mathbb{E}({\mathcal R}_{1,t}^n)\right]> r\right)&=&\mathbb{P}\left(\exp\left(\lambda n\sup_{t\in [0,1]}\left[{\mathcal R}_{1,t}^n-\mathbb{E}({\mathcal R}_{1,t}^n)\right]\right)> e^{n\lambda r}\right)\\
&\le&e^{-n\lambda r}\mathbb{E}\left(\exp\left(\lambda n\left[{\mathcal R}_{1,1}^n-\mathbb{E}({\mathcal R}_{1,1}^n)\right]\right)\right)
\end{eqnarray*}
and similarly
$$\mathbb{P}\left(\inf_{t\in [0,1]}\left[{\mathcal R}_{q,t}^n-\mathbb{E}({\mathcal R}_{1,t}^n)\right]<-r\right)\le e^{-n\lambda r}\mathbb{E}\left(\exp\left(-\lambda n\left[{\mathcal R}_{1,1}^n-\mathbb{E}({\mathcal R}_{1,1}^n)\right]\right)\right).
$$

So we get
\begin{equation*}\dfrac{1}{n} \log \mathbb{P}\left(\sup_{t\in [0,1]}\left[{\mathcal R}_{1,t}^n-\mathbb{E}({\mathcal R}_{1,t}^n)\right]> r \right)\le -\lambda r + \frac 1n\log\mathbb{E}\left(e^{\lambda n {\mathcal R}_{1,1}^n}\right)-\lambda\mathbb{E}({\mathcal R}_{1,1}^n).
\end{equation*}

It is easy to see that $\mathbb{E}({\mathcal R}_{1,1}^n) \rightarrow 0$ as $n$ goes to infinity. We have already seen in (\ref{K}) that $\frac 1n\log\mathbb{E}\left(e^{\lambda n {\mathcal R}_{1,1}^n}\right)\rightarrow 0$ as $n$ gos to infinity. So we obtain for all $\lambda>0$
\begin{equation*}\limsup_{n\rightarrow \infty}\dfrac{1}{n} \log \mathbb{P}\left(\sup_{t\in [0,1]}\left[{\mathcal R}_{1,t}^n-\mathbb{E}({\mathcal R}_{1,t}^n)\right]> r \right)\le -\lambda r.
\end{equation*}
Letting $\lambda>0$ goes to infinity, we obtain that the left term in the last inequality goes to $-\infty$.

And similarly, by doing the same calculations with 
$$\mathbb{P}\left(\inf_{t\in [0,1]}\left[{\mathcal R}_{q,t}^n-\mathbb{E}({\mathcal R}_{1,t}^n)\right]<-r\right),
$$
we obtain that
\begin{equation*}
\sup_{t\in [0,1]}\left|{\mathcal R}_{1,t}^n-\mathbb{E}({\mathcal R}_{1,t}^n)\right|\superexpldp 0.
\end{equation*}

Since
$$\mathbb{E}({\mathcal R}_{1,1}^n) \superexpldp 0,$$
we obtain (\ref{negR123}) for $q=1$. 
\vspace{10pt}

\subsection{Proof of Theorem  \ref{MDP1}}$\quad$
\vspace{10pt}

As is usual, the proof of the MDP is somewhat simpler than the LDP, relying on the same line of proof than the one for the CLT. Namely, a good control of the asymptotic of the moment generating functions, and G\"artner-Ellis theorem. We shall then prove that for all $\lambda\in \mathbb{R}^3$
\begin{equation}
\label{MDP1-prove}
 \lim_{n \to \infty}\dfrac{1}{b_n^2}\log\mathbb{E}\exp\left(b_n^2\dfrac{\sqrt{n}}{b_n} \left\langle \lambda,V_1^n - [V]_1 \right\rangle\right)=\dfrac{1}{2}\left\langle\lambda,\Sigma_{1}\cdot\lambda\right\rangle.
\end{equation}

Taking the calculation in \eqref{lem_2},we have
$$\dfrac{1}{b_n^2}\log\mathbb{E}\exp\left(b_n^2\dfrac{\sqrt{n}}{b_n} \left\langle \lambda,V_1^n - [V]_1 \right\rangle\right)=\dfrac{1}{b_n^2}\sum_{k=1}^{n}\left[H_k^n(\lambda)-b_{n}\sqrt{n}\left\langle \lambda,[V]_1 \right\rangle\right],$$
with
$$H_k^n(\lambda):=P_{c_k^n}\left(\lambda_{1}b_{n}\sqrt{n}a_{1,k},\lambda_{2}b_{n}\sqrt{n}a_{2,k},\lambda_{3}b_{n}\sqrt{n}\sqrt{a_{1,k}}\sqrt{a_{2,k}}\right),$$
where $a_{\ell,k}$ are given in (\ref{defxi}) and $c_k^n$ is given in (\ref{defc}).

By our condition \eqref{Max_condition}, 
$$\varepsilon(n):=\sqrt{n}b_n\max_{1 \leqslant k \leqslant n}\max_{\ell =1,2}a_{\ell,k}\xrightarrow[{n \to \infty}]\quad 0.$$

By multidimensional Taylor formula and noting that $P_{c_k^n}(0,0,0)=0$, $\nabla P_{c_k^n}(0,0,0)=(1,1,c_k^n)^T$ and the Hessian matrix
$$
 H(P_{c_k^n})(0,0,0)= \begin{pmatrix}
                                       \vspace{0.3cm}
                                       2 & 2(c_k^n)^2 & 2c_k^n\\
                                       \vspace{0.3cm}
                                       2(c_k^n)^2 & 2 & 2c_k^n\\
                                       2c_k^n & 2c_k^n & 1+(c_k^n)^2
                                       \end{pmatrix},
 $$
and after an easy calculations, we obtain once if $||\lambda||\cdot|\varepsilon(n)|<\frac 14,$ i.e. for $n$ large enough,
\begin{eqnarray*}
H_k^n(\lambda)=b_{n}\sqrt{n}\left\langle \lambda,[V]_1 \right\rangle +n b_n^2\frac 12\langle \lambda, \Sigma_k^n\cdot \lambda  \rangle+n b_n^2\nu(k,n),
\end{eqnarray*}
where
$$
 \Sigma_k^n:= \begin{pmatrix}
                                       \vspace{0.3cm}
                                       a_{1,k}^2 & (\vartheta_k^n)^2  & a_{1,k}\vartheta_k^n  \\
                                       \vspace{0.3cm}
                                       (\vartheta_k^n)^2  & a_{2,k}^{2}  & a_{2,k}\vartheta_k^n \\
                                       a_{1,k}\vartheta_k^n & a_{2,k}\vartheta_k^n   & \dfrac{1}{2}(a_{1,k}a_{2,k}+(\vartheta_k^n)^2)
                                       \end{pmatrix},
 $$
where $\vartheta_k^n$ is given in (\ref{defc}), and $\nu(k,n)$ satisfies
 $$|\nu(k,n)| \leqslant C ||\lambda||\cdot|\varepsilon(n)|,$$
where $C=\frac{1}{6}\sup_{||\lambda||\le 1/4}\biggl|\dfrac{\partial^3 P(\lambda_1,\lambda_2,\lambda_3)}{\partial^3\lambda}\biggr|$.

On the other hand, by the classical Lebesgue derivation theorem see \cite{Djellout1}, we have
$$
\sum_{k=1}^{n}\left(\dfrac{\int_{t_{k-1}^n}^{t_k^n}g(s) \mathrm ds}{t_{k}^{n}-t_{k-1}^{n}}\right)\left(\dfrac{\int_{t_{k-1}^n}^{t_k^n}h(s) \mathrm ds}{t_{k}^{n}-t_{k-1}^{n}}\right)(t_{k}^{n}-t_{k-1}^{n}) \rightarrow \int_{0}^{1}g(s)h(s)\mathrm ds
$$
by taking different chosse of $g$ and $h$: once $g(s)=h(s)=\sigma_{\ell,s}^2$ or $g(s)=h(s)=\sigma_{1,s}\sigma_{1,s}\rho_s$, or  $g(s)=\sigma_{\ell,s}^2$ and $h(s)=\sigma_{\ell',s}^2$ $\ell\not=\ell'$, and  $g(s)=\sigma_{\ell,s}^2$ and $h(s)=\sigma_{1,s}\sigma_{1,s}\rho_s$, we obtain that

$$n\sum_{k=1}^{n}\Sigma_k^n \rightarrow_{n \to \infty}
 \Sigma_{1}= \begin{pmatrix}
                                       \vspace{0.3cm}
                                       \int_0^1\sigma_{1,t}^4 \mathrm dt &  \int_0^1\sigma_{1,t}^2\sigma_{2,t}^2\rho_t^2 \mathrm dt & \int_0^1\sigma_{1,t}^3\sigma_{2,t}\rho_t \mathrm dt\\
                                        \vspace{0.3cm}
                                        \int_0^1\sigma_{1,t}^2\sigma_{2,t}^2\rho_t^2 \mathrm dt &  \int_0^1\sigma_{2,t}^4 \mathrm dt & \int_0^1\sigma_{1,t}\sigma_{2,t}^3\rho_t \mathrm dt\\
                                        \int_0^1\sigma_{1,t}^3\sigma_{2,t}\rho_t \mathrm dt & \int_0^1\sigma_{1,t}\sigma_{2,t}^3\rho_t \mathrm dt &  \int_0^1\dfrac{1}{2}\sigma_{1,t}^2\sigma_{2,t}^2(1+\rho_t^2) \mathrm dt
                                      \end{pmatrix}=\Sigma_1.
$$ 

Then the \eqref{MDP1-prove} follows.

Hence \eqref{MDP1} follows from the G$\ddot a$rtner-Ellis theorem.

\vspace{10pt}

\subsection{Proof of Theorem  \ref{MDP2}}$\quad$
\vspace{10pt}

It is well known that the LDP of finite dimensional vector
$$\left(\dfrac{\sqrt{n}}{b_n}(V_{s_1}^n(X)-[V]_{s_1},\cdots,V_{s_k}^n(X)-[V]_{s_k}) \right), 0< s_1 < \cdots < s_k \leqslant 1,k\geqslant 1 $$
and the following exponential tightness: for any $s\in [0, 1]$ and $\eta > 0$
$$\lim_{\varepsilon \downarrow 0}\limsup_{n \to \infty}\dfrac{1}{b_n^2} \log \mathbb{P}\left(\dfrac{\sqrt{n}}{b_n}\sup_{s \leqslant t \leqslant s+\varepsilon}\|\Delta_s^t V_{\cdot}^n(X) - \Delta_s^t [V]_{\cdot}\| > \eta \right)=-\infty$$
with $\Delta_s^t V_{\cdot}^n=V_t^n - V_s^n$, are sufficient for the LDP of $\dfrac{\sqrt{n}}{b_n}(V_{\cdot}^n(X)-[V]_{\cdot})$ for the sup-norm topology (cf. \cite{Demzei},\cite{Djellout1}).

\vspace{10pt}
Under the assumption of Theorem \ref{MDP2}, we have:
\begin{equation}\label{E}\lim_{n \to \infty} \sqrt{n}\sup_{t \in [0,1]}\|\mathbb{E}V_t^n(X)-[V]_t\|=0.\end{equation}

In fact, we have:
\begin{align*}
&\sqrt{n}\sup_{t \in [0,1]}\|\mathbb{E}V_t^n(X)-[V]_t\| \\
&\leqslant 3\sqrt{n}\max\left(\max_{\ell=1,2}\sup_{t \in [0,1]}|\mathbb{E}Q_{\ell,t}^n(X)-[X_{\ell}]_t|,\sqrt{n}\sup_{t \in [0,1]}|\mathbb{E}C_{t}^n(X)-\langle X_1,X_2\rangle_t|\right)\\
&\leqslant 3\max\left(\max_{\ell=1,2}\max_{k \leqslant n}\sqrt{\int_{(k-1)/n}^{k/n}\sigma_{\ell,t}^4\mathrm dt},\max_{k \leqslant n}\sqrt{\int_{(k-1)/n}^{k/n}\sigma_{1,t}^2\sigma_{2,t}^2\rho_t^2\mathrm dt}\right).
\end{align*}

By our condition \eqref{Max_condition}, we obtain (\ref{E}).
\vspace{10pt}

Now, we show that for any partition $0< s_1 < \cdots < s_k \leqslant 1,k\geqslant 1$ of $[0,1]$
$$\dfrac{\sqrt{n}}{b_n}\left(V_{s_1}^n(X)-[V]_{s_1},\cdots,V_{s_k}^n(X)-[V]_{s_k}\right)$$
satisfies the LDP on $\mathbb R^k$ with speed $b_n^2$ and with the rate function given by
\begin{equation}
\label{LDP_multi}
  I_{s_1,\cdots,s_k}(x_1,\cdots,x_k)=\dfrac{1}{2}\sum_{i=1}^{k}\langle(x_i-x_{i-1}),(\Sigma_{s_{i-1}}^{s_i})^{-1}\cdot (x_i-x_{i-1})\rangle 
\end{equation}
where
$$
 \Sigma_s^u= \begin{pmatrix}
                                       \vspace{0.3cm}
                                       \int_s^u\sigma_{1,t}^4 \mathrm dt &  \int_s^u\sigma_{1,t}^2\sigma_{2,t}^2\rho_t^2 \mathrm dt & \int_s^u\sigma_{1,t}^3\sigma_{2,t}\rho_t \mathrm dt\\
                                        \vspace{0.3cm}
                                        \int_s^u\sigma_{1,t}^2\sigma_{2,t}^2\rho_t^2 \mathrm dt &  \int_s^u\sigma_{2,t}^4 \mathrm dt & \int_s^u\sigma_{1,t}\sigma_{2,t}^3\rho_t \mathrm dt\\
                                        \int_s^u\sigma_{1,t}^3\sigma_{2,t}\rho_t \mathrm dt & \int_s^u\sigma_{1,t}\sigma_{2,t}^3\rho_t \mathrm dt &  \int_s^u\dfrac{1}{2}\sigma_{1,t}^2\sigma_{2,t}^2(1+\rho_t^2) \mathrm dt
                                      \end{pmatrix}
 $$
is invertible and 
\begin{align*}
{\rm det}(Q_s^u) &= \biggl(\int_s^u\sigma_{1,t}^4 \mathrm dt\biggr)\biggl(\int_s^u\sigma_{2,t}^4 \mathrm dt\biggr)\biggl(\int_s^u\dfrac{1}{2}\sigma_{1,t}^2\sigma_{2,t}^2(1+\rho_t^2) \mathrm dt\biggr)\\
  &+2\biggl(\int_s^u\sigma_{1,t}^2\sigma_{2,t}^2\rho_t^2 \mathrm dt\biggr)\biggl(\int_s^u\sigma_{1,t}\sigma_{2,t}^3\rho_t \mathrm dt\biggr)\biggl(\int_s^u\sigma_{1,t}^3\sigma_{2,t}\rho_t \mathrm dt\biggr)\\
  &-\biggl(\int_s^u\dfrac{1}{2}\sigma_{1,t}^2\sigma_{2,t}^2(1+\rho_t^2) \mathrm dt\biggr)\biggl(\int_s^u\sigma_{1,t}^2\sigma_{2,t}^2\rho_t^2 \mathrm dt\biggr)^2\\
  &-\biggl(\int_s^u\sigma_{1,t}^4 \mathrm dt\biggr)\biggl(\int_s^u\sigma_{1,t}\sigma_{2,t}^3\rho_t \mathrm dt\biggr)^2-\biggl(\int_s^u\sigma_{2,t}^4 \mathrm dt\biggr)\biggl(\int_s^u\sigma_{1,t}^3\sigma_{2,t}\rho_t \mathrm dt\biggr)^2
 \end{align*}  
and $(\Sigma_s^u)^{-1}$ his inverse.

For $n$ large enough we have $1<[nt_1]<\cdots<[nt_k]< n$, so by applying the homeomorphism 
$$\Upsilon : (x_1,\cdots,x_k) \to (x_1,x_2-x_1,\cdots,x_k-x_{k-1})$$
$Z_n=(V_{s_1}^n(X)-[V]_{s_1},\cdots,V_{s_k}^n(X)-[V]_{s_k})$ can be mapped to $U_n=\Upsilon Z_n$ with independent components. 

Then we consider the LDP of  $\dfrac{\sqrt n}{b_n}(U_n-\mathbb{E}U_n)$.

For any $\theta = (\theta_1,\cdots,\theta_k)\in (\mathbb{R}^3)^k$,
$$\Lambda{s_1,\cdots,s_k}(\theta)=\lim_{n \to \infty}\dfrac{1}{b_n^2}\log\mathbb{E}\exp\left(b_n\sqrt{n}\left\langle \lambda,U_n-EU_n \right\rangle \right)=\sum_{i=1}^{k}\frac{1}{2}\langle \lambda_i,\Sigma_{s_{i-1}}^{s_i}\cdot\lambda_i\rangle.$$

By G$\ddot a$rtner-Ellis theorem, $\dfrac{\sqrt n}{b_n}(U_n-\mathbb{E}U_n)$ satifies the LDP in $(\mathbb{R}^3)^k$ with speed $b_n^2$ and with the good rate function
$$\Lambda^*{s_1,\cdots,s_k}(x)=\dfrac{1}{2}\sum_{i=1}^{k}\langle x_i,(\Sigma_{s_{i-1}}^{s_i})^{-1}\cdot x_i\rangle.
$$

Then by the inverse contraction principle, we have $\dfrac{\sqrt n}{b_n}(Z_n-\mathbb{E} Z_n)$ satisfies the LDP with speed $b_n^2$ and with the rate function $I_{s_1,\cdots,s_k}(x)$ given in (\ref{LDP_multi}).

\vspace{10pt}

Now, we shall prove that for any $\eta>0$, $s\in [0,1]$
\begin{equation} 
\label{limsup}
\lim_{\varepsilon \downarrow 0}\limsup_{n \to \infty}\dfrac{1}{b_n^2} \log \mathbb{P}\left(\dfrac{\sqrt{n}}{b_n}\sup_{s \leqslant t \leqslant s+\varepsilon}\|\Delta_s^t V_{\cdot}^n(X) - \mathbb{E}\Delta_s^t V_{\cdot}^n(X)\| > \eta \right)=-\infty.
\end{equation}

For that we need to prove that for $\ell=1,2$ and for all $\eta>0$ and $s\in [0,1]$
\begin{equation}
\label{yac1}
\lim_{\varepsilon \downarrow 0}\limsup_{n \to \infty}\dfrac{1}{b_n^2} \log \mathbb{P}\left(\dfrac{\sqrt{n}}{b_n}\sup_{s \leqslant t \leqslant s+\varepsilon}|\Delta_s^t Q_{\ell,\cdot}(X) - \mathbb{E}\Delta_s^t Q_{\ell,\cdot}(X)| > \eta \right)=-\infty,
\end{equation}
and
\begin{equation}
\label{yac0}
\lim_{\varepsilon \downarrow 0}\limsup_{n \to \infty}\dfrac{1}{b_n^2} \log \mathbb{P}\left(\dfrac{\sqrt{n}}{b_n}\sup_{s \leqslant t \leqslant s+\varepsilon}|\Delta_s^t C_{\cdot}(X) - \mathbb{E}\Delta_s^t C_{\cdot}(X)| > \eta\right)=-\infty.
\end{equation}

In fact  (\ref{yac1}) can be done in the same way than in Djellout et al.\cite{Djellout1}. It remains to show (\ref{yac0}). This will be done following the same technique as for the proof of (\ref{yac1}) and using a result of \cite{Djellout2}. Remark that $(C^n_{t}(X) - \mathbb{E}C^n_{t}(X))$ is an $\mathcal{F}_{[nt]/n}$-martingale. Then  
$$\exp(\lambda[\Delta_s^t (C^n_{\cdot}(X) - \mathbb{E}C^n_{\cdot}(X))])
$$
is a sub-martingale. By the maximal inequality, we have for any $\eta, \lambda > 0$
\begin{eqnarray}\label{yac2}
\mathbb{P}\left(\sup_{s \leqslant t \leqslant s+\varepsilon}\Delta_s^t \left[C^n_{\cdot}(X) -\mathbb{E}C^n_{\cdot}(X)\right] > \eta\right) &=& \mathbb{P}\left(\exp(\lambda\sup_{s \leqslant t \leqslant s+\varepsilon}\Delta_s^t \left[C^n_{\cdot}(X) - \mathbb{E}C^n_{\cdot}(X)\right] > e^{\lambda\eta} \right)\nonumber\\
 &\leqslant&  e^{-\lambda\eta}\mathbb{E}\exp\left(\lambda \Delta_s^{s+\varepsilon} \left[C^n_{\cdot}(X) - \mathbb{E}C^n_{\cdot}(X)\right]\right),
\end{eqnarray}
and similary,
\begin{equation}\label{yac3}
\mathbb{P}\left(\inf_{s \leqslant t \leqslant s+\varepsilon}\Delta_s^t \left[C^n_{\cdot}(X) - \mathbb{E}C^n_{\cdot}(X)\right] < -\eta\right) \leqslant  e^{-\lambda\eta}\mathbb{E}\exp\left(-\lambda[\Delta_s^{s+\varepsilon} \left[C^n_{\cdot}(X) - \mathbb{E}C^n_{\cdot}(X)\right]\right).
\end{equation}

Using Remark 2.4 in \cite{Djellout2}, we have that for all $c \in \mathbb{R}$
$$\lim_{n\to \infty}\dfrac{1}{b_n^2}\log \mathbb{E}\exp\left(cb_n^2\dfrac{\sqrt n}{b_n}\Delta_s^{s+\epsilon}\left[C^n_{\cdot}(X) - \mathbb{E}C^n_{\cdot}(X)\right]\right)=\dfrac{1}{2}c^2 \int_s^{s+\epsilon} \sigma_{1,t}^2\sigma_{2,t}^2(1+\rho_t^2) \mathrm dt.$$

Therefore taking $\eta = \delta\dfrac{b_n}{n}$, $\lambda=b_n\sqrt{n}c$ $(c>0)$ in~\eqref{yac2}, we get
\begin{eqnarray*}
&&\limsup_{n \to \infty}\dfrac{1}{b_n^2} \log \mathbb{P}\left(\dfrac{\sqrt{n}}{b_n}\sup_{s \leqslant t \leqslant s+\varepsilon}\Delta_s^t \left[C^n_{\cdot}(X) - \mathbb{E}C^n_{\cdot}(X)\right] > \delta \right)\\
&&\leqslant \inf_{c>0}\left\{-c\delta + \dfrac{1}{2}c^2 \int_s^{s+\epsilon} \sigma_{1,t}^2\sigma_{2,t}^2(1+\rho_t^2) \mathrm dt \right\}=- \dfrac{\delta^2}{2 \int_s^{s+\epsilon} \sigma_{1,t}^2\sigma_{2,t}^2(1+\rho_t^2) \mathrm dt},
\end{eqnarray*}
and similary by \eqref{yac3},
$$
\limsup_{n \to \infty}\dfrac{1}{b_n^2} \log \mathbb{P}\left(\dfrac{\sqrt{n}}{b_n}\inf_{s \leqslant t \leqslant s+\varepsilon}\Delta_s^t \left[C^n_{\cdot}(X) - \mathbb{E}C^n_{\cdot}(X)\right]<-\delta \right)\leqslant - \dfrac{\delta^2}{2 \int_s^{s+\epsilon} \sigma_{1,t}^2\sigma_{2,t}^2(1+\rho_t^2) \mathrm dt}.
$$

By the integrability of $\sigma_{1,t}^2\sigma_{2,t}^2(1+\rho_t^2)$, we have
$$
\lim_{\varepsilon \downarrow 0}\sup_{s\in [0,1]}\int_s^{s+\epsilon}\sigma_{1,t}^2\sigma_{2,t}^2(1+\rho_t^2) \mathrm dt=0.
$$

Hence \eqref{yac0} follows from the above estimations. So we have \eqref{limsup}.
\vspace{10pt}

By \eqref{LDP_multi} and \eqref{limsup}, $\dfrac{\sqrt{n}}{b_n}(V_{\cdot}^n-[V]_{\cdot})$ 
satifies the LDP with the speed $b_n^2$ and the good rate function
 $$
I_{sup}(x)=\sup \left\{I_{s_1,\cdots,s_k}(x(s_1),\cdots,x(s_k));  0< s_1 < \cdots < s_k \leqslant 1,k\geqslant 1\right\},
$$
where
$$
I_{s_1,\cdots,s_k}(x(s_1),\cdots,x(s_k)) = \dfrac{1}{2}\sum_{i=1}^{k}\langle x(s_i)-x(s_{i-1}),(\Sigma_{s_{i-1}}^{s_i})^{-1}\cdot (x(s_i)-x(s_{i-1}))\rangle
$$
\vspace{5pt}

It remains to prove that $I_{sup}(x)=I_{mdp}(x)$.

We shall prove that $I_{sup}(x)\leqslant I_{mdp}(x).$

For this, we treat the first element of the matrix $(\Sigma_{s_{i-1}}^{s_i})^{-1}$ which is denoted $(\Sigma_{s_{i-1}}^{s_i})^{-1}_{1,1}$ and we prove that
$$\sum_{i=1}^{k}(x_1(s_i)-x_1(s_{i-1}))^2.(\Sigma_{s_{i-1}}^{s_i})^{-1}_{1,1}\leqslant \int_{0}^{1}(x_1'(t))^2.(\Sigma_{t}^{-1})_{1,1}\mathrm dt,
$$
where $(\Sigma_{t}^{-1})_{1,1}$ represente the first element of the matrix $\Sigma_{t}^{-1}$.

We have
$$(\Sigma_{s_{i-1}}^{s_i})^{-1}_{1,1}=\frac{1}{{\rm det}(\Sigma_{s_{i-1}}^{s_i})}\biggl(\biggl(\int_{s_{i-1}}^{s_i}\sigma_{2,t}^4 \mathrm dt\biggr)\biggl(\int_{s_{i-1}}^{s_i}\dfrac{1}{2}\sigma_{1,t}^2\sigma_{2,t}^2(1+\rho_t^2) \mathrm dt\biggr)-\biggl(\int_{s_{i-1}}^{s_i}\sigma_{1,t}\sigma_{2,t}^3\rho_t \mathrm dt\biggr)^2\biggr).$$

By \cite[p.1305]{Jiang1}, for $x:=(x_1,x_2,x_3) \in \mathcal H$, if $I_{sup}(x)< +\infty$, then for $0< s_1 < \cdots < s_k \leqslant 1$,

Then
\begin{align*}
 \sum_{i=1}^{k}(x_1(s_i)-x_1(s_{i-1}))^2.(\Sigma_{s_{i-1}}^{s_i})^{-1}_{1,1}&\leqslant \int_{0}^{1}(x'_1(t))^2.\frac{\dfrac{1}{2}\sigma_{1,t}^2\sigma_{2,t}^6(1+\rho_t^2)-\sigma_{1,t}^2\sigma_{2,t}^6\rho_t^2}{{\rm det}(\Sigma_{t})} \mathrm dt\\
 &= \int_{0}^{1}(x_1'(t))^2.(\Sigma_{t}^{-1})_{1,1} \mathrm dt,
\end{align*}

The same calculation with the other terms of the matrix given in the following, implies that $I_{sup}(x)\leqslant I_{mdp}(x)$:

$$(\Sigma_{s_{i-1}}^{s_i})^{-1}_{2,2}=\frac{1}{{\rm det}(\Sigma_{s_{i-1}}^{s_i})}\biggl[\biggl(\int_{s_{i-1}}^{s_i}\sigma_{1,t}^4 \mathrm dt\biggr)\biggl(\int_{s_{i-1}}^{s_i}\dfrac{1}{2}\sigma_{1,t}^2\sigma_{2,t}^2(1+\rho_t^2) \mathrm dt\biggr)-\biggl(\int_{s_{i-1}}^{s_i}\sigma_{1,t}^3\sigma_{2,t}\rho_t \mathrm dt\biggr)^2\biggr],$$

$$(\Sigma_{s_{i-1}}^{s_i})^{-1}_{3,3}=\frac{1}{{\rm det}(\Sigma_{s_{i-1}}^{s_i})}\biggl[\biggl(\int_{s_{i-1}}^{s_i}\sigma_{1,t}^4 \mathrm dt\biggr)\biggl(\int_{s_{i-1}}^{s_i}\sigma_{2,t}^4 \mathrm dt\biggr)-\biggl(\int_{s_{i-1}}^{s_i}\sigma_{1,t}^2\sigma_{2,t}^2\rho_t^2 \mathrm dt\biggr)^2\biggr],$$

\begin{eqnarray*}
&&(\Sigma_{s_{i-1}}^{s_i})^{-1}_{1,2}=(\Sigma_{s_{i-1}}^{s_i})^{-1}_{2,1}=\frac{1}{{\rm det}(\Sigma_{s_{i-1}}^{s_i})}\biggl[\biggl(\int_{s_{i-1}}^{s_i}\sigma_{1,t}^3\sigma_{2,t}\rho_t \mathrm dt\biggr)\biggl(\int_{s_{i-1}}^{s_i}\sigma_{1,t}\sigma_{2,t}^3\rho_t \mathrm dt\biggr)\\
&&\qquad \qquad\qquad\qquad\qquad-\biggl(\int_{s_{i-1}}^{s_i}\sigma_{1,t}^2\sigma_{2,t}^2\rho_t^2 \mathrm dt\biggr)\biggl(\int_{s_{i-1}}^{s_i}\dfrac{1}{2}\sigma_{1,t}^2\sigma_{2,t}^2(1+\rho_t^2) \mathrm dt\biggr)\biggr],
\end{eqnarray*}

\begin{eqnarray*}
&&(\Sigma_{s_{i-1}}^{s_i})^{-1}_{1,3}=(\Sigma_{s_{i-1}}^{s_i})^{-1}_{3,1}=\frac{1}{{\rm det}(\Sigma_{s_{i-1}}^{s_i})}\biggl[\biggl(\int_{s_{i-1}}^{s_i}\sigma_{1,t}^2\sigma_{2,t}^2\rho_t^2 \mathrm dt\biggr)\biggl(\int_{s_{i-1}}^{s_i}\sigma_{1,t}\sigma_{2,t}^3\rho_t \mathrm dt\biggr)\\
&&\qquad \qquad\qquad\qquad\qquad-\biggl(\int_{s_{i-1}}^{s_i}\sigma_{2,t}^4 \mathrm dt\biggr)\biggl(\int_{s_{i-1}}^{s_i}\sigma_{1,t}^3\sigma_{2,t}\rho_t \mathrm dt\biggr)\biggr],
\end{eqnarray*}

\begin{eqnarray*}
&&(\Sigma_{s_{i-1}}^{s_i})^{-1}_{2,3}=(\Sigma_{s_{i-1}}^{s_i})^{-1}_{3,2}=\frac{1}{{\rm det}(\Sigma_{s_{i-1}}^{s_i})}\biggl[\biggl(\int_{s_{i-1}}^{s_i}\sigma_{1,t}^2\sigma_{2,t}^2\rho_t^2 \mathrm dt\biggr)\biggl(\int_{s_{i-1}}^{s_i}\sigma_{1,t}^3\sigma_{2,t}\rho_t \mathrm dt\biggr)\\
&&\qquad \qquad\qquad\qquad\qquad-\biggl(\int_{s_{i-1}}^{s_i}\sigma_{1,t}^4 \mathrm dt\biggr)\biggl(\int_{s_{i-1}}^{s_i}\sigma_{1,t}\sigma_{2,t}^3\rho_t \mathrm dt\biggr)\biggr].
\end{eqnarray*}

On the other hand, by the convergence of martingales and Fatou's lemma,
$$I_{mdp}(x)< +\infty, \qquad {\rm and}\qquad I_{mdp}(x)\leqslant I_{sup}(x).$$

So we have $I_{sup}(x)=I_{mdp}(x).$

\subsection{Proof of  Theorem \ref{Lipschitz} } $\qquad$

\vspace{10pt}

\underline{ Step 1}. We shall prove that $\tilde V_{\cdot}^n$ and $V_{\cdot}^n(X-Y)$ satisfy the same LDP, by means of the approximation Lemma in \cite{Djellout1}. So we shall prove that  $\tilde Q_{\ell,\cdot}^n$ and $Q_{\ell,\cdot}^n(X-Y)$ satisfy the same LDP and idem for $\tilde C_{\cdot}^n$ and $C_{\cdot}^n(X-Y)$. 

We have
\begin{equation}
\label{Drift1}
\sup_{t\in [0,1]}|\tilde Q_{\ell,t}^n(X) - Q_{\ell,t}^n(X-Y)|\leqslant \varepsilon(n)Q_{\ell,t}^n(X-Y)+\left(1+\frac{1}{\varepsilon(n)}\right)Z_{\ell,n}
\end{equation}
and
\begin{equation}
\label{Drift3}
\sup_{t\in [0,1]}|\tilde C_{1,t}^n(X) - C_{1,t}^n(X-Y)| \leqslant \varepsilon(n)\sum_{\ell=1}^2Q_{\ell,t}^n(X-Y)+\left(\frac 12+\frac{1}{\varepsilon(n)}\right)\sum_{\ell=1}^2Z_{\ell,n},
\end{equation}
where the sequence $\varepsilon(n)>0$ will be selected later, and $Z_{\ell,n}$ is given 
\begin{equation}\label{defZ}
Z_{\ell,n}=\sum_{k=1}^{n}\left(\int_{t_{k-1}^n}^{t_k^n}b_{\ell,t}(X_{t})\mathrm dt - b_{\ell,t_{k-1}^n}(X_{t_{k-1}^n})(t_{k}^n-t_{k-1}^n) \right)^2,
\end{equation}
with $X_t=(X_{1,t}, X_{2,t})$.

For $Q_{\ell,t}^n(X-Y)$, being a Gaussian process, Theorem 1.1 in \cite{Djellout1} may be used. It remains to control $Z_{\ell,n}$. For this we just need to prove that:
\begin{equation} 
\label{Drift_1}
\dfrac{1}{\varepsilon(n)}Z_{\ell,n}\superexpldp 0.
\end{equation}

The main idea is to reduce it to estimations of $M_{\ell,t}=\int_0^t\sigma_{\ell,s}\mathrm dB_{\ell,s}$, by
means of Gronwall's inequality. So, we have at first for all $t \in [0,1]$
\begin{align*}
\|X_t\| &\leqslant \|X_0\|+C \int_0^t (1+(1+\eta(s))\|X_s\|) \mathrm ds + \sup_{s \leqslant t}\|M_s\| \\
&\leqslant \left(C+\|X_{1,0}\|+\|X_{2,0}\|+\sup_{s \leqslant 1}\|M_s\| \right) + C_1 \int_0^t \|X_s\| \mathrm ds.
\end{align*}
where $C_1=C(1+\eta(1))$. Hence, by Gronwall's inequality
\begin{equation} 
\label{Gronw_1}
\|X_t\| \leqslant \left(C+\|X_{1,0}\|+\|X_{2,0}\|+\sup_{s \leqslant 1}\|M_s\| \right)e^{C_1t}, \qquad \forall t \in [0,1]
\end{equation}

For any $s \in [0,1]$, $v>0$
\begin{eqnarray} 
\label{Gronw_2}
\sup_{s \leqslant t \leqslant s+v}\|X_t-X_s\| &\leqslant&\sup_{s \leqslant t \leqslant s+v}\|M_t-M_s\|+v.\sup_{s \leqslant t \leqslant s+v}\|b(t,X_t)\| \nonumber\\
&\leqslant&\sup_{s \leqslant t \leqslant s+v}\|M_t-M_s\|+vC_2\left(\sup_{0 \leqslant t \leqslant 1}\|X_t\|+1\right)
\end{eqnarray}

We get by \eqref{yac_drift}, \eqref{Gronw_1},\eqref{Gronw_2} and Cauchy-Schwarz's inequality

\begin{eqnarray} 
\label{Drift_4}
&&\left(\int_{t_{k-1}^n}^{t_k^n}b_{\ell}(t,X_{t})\mathrm dt - b_{\ell}(t_{k-1}^n, X_{t_{k-1}^n})(t_{k}^n-t_{k-1}^n) \right)^2\nonumber\\
&\leqslant& \left(\dfrac{1}{n}C\left(1+ \sup_{t_{k-1}^n \leqslant t \leqslant t_k^n}\|X_t-X_{t_{k-1}^n}\|+2\eta\left(\frac{1}{n}\right)\sup_{0 \leqslant t \leqslant 1}\|X_t\| \right)\right)^2\nonumber\\
&\leqslant& \dfrac{C_3}{n^2}\left(1+ \sup_{t_{k-1}^n \leqslant t \leqslant t_k^n}\|M_t-M_{t_{k-1}^n}\|^2+\left(\frac{1}{n^2}+\eta\left(\frac{1}{n}\right)^2\right)\sup_{0 \leqslant t \leqslant 1}\|M_t\|^2\right)
\end{eqnarray}

 Chose $\varepsilon(n)>0$ so that 
 \begin{equation} 
\label{eps}
\varepsilon(n) \to 0 \quad {\rm but }\quad \frac{\frac{1}{n^2}+\eta\left(\frac{1}{n}\right)^2}{\varepsilon(n)} \to 0
\end{equation}

By \eqref{Drift_4} and the definition of $Z_{\ell,n}$, we have that
\begin{equation*} 
\limsup_{n \to \infty}\dfrac{1}{n} \log \mathbb{P}\left(\dfrac{1}{\varepsilon(n)}Z_{\ell,n}>\delta \right)\leqslant \max(A,B) 
\end{equation*}
where
\begin{eqnarray} 
\label{A_1}
A &=&\limsup_{n \to \infty}\dfrac{1}{n} \log \mathbb{P}\left(\dfrac{1}{\varepsilon(n)n}\max_{k\leqslant n}\sup_{t_{k-1}^n \leqslant t \leqslant t_k^n}\|M_t-M_{t_{k-1}^n}\|^2>C_4\delta \right)\nonumber\\
&\leqslant& 2\max_{\ell=1,2} \limsup_{n \to \infty}\dfrac{1}{n} \log \mathbb{P}\left(\dfrac{1}{\varepsilon(n)n}\max_{k\leqslant n}\sup_{t_{k-1}^n \leqslant t \leqslant t_k^n}|M_{\ell,t}-M_{\ell,t_{k-1}^n}|^2>C_4\delta \right),
\end{eqnarray}
and
\begin{eqnarray} 
\label{B_1}
B &=&\limsup_{n \to \infty}\dfrac{1}{n} \log \mathbb{P}\left(\dfrac{1}{\varepsilon(n)^2n}\left(\frac{1}{n^2}+\eta\left(\frac{1}{n}\right)^2\right)\sup_{0 \leqslant t \leqslant 1}\|M_t\|^2>C_5\delta \right)\nonumber\\
&\leqslant& 2\max_{\ell=1,2} \limsup_{n \to \infty}\dfrac{1}{n} \log \mathbb{P}\left(\dfrac{1}{\varepsilon(n)n}\left(\frac{1}{n^2}+\eta\left(\frac{1}{n}\right)^2\right)\sup_{0 \leqslant t \leqslant 1}|M_{\ell,t}|^2>C_5\delta \right).
\end{eqnarray}

By L\'evy's inequality for a Brownian motion and our choice \eqref{eps} of $\varepsilon(n)$, the limits \eqref{A_1} and \eqref{B_1} are both $-\infty$. Limit (\ref{Drift_1}) follows.

\vspace{10pt}

\underline{Step 2}. We shall prove that $\frac{\sqrt n}{b_n}(\tilde V_{\cdot}^n-[V]_{\cdot})$ and $\frac{\sqrt n}{b_n}(V_{\cdot}^n(X-Y)-[V]_{\cdot})$ satisfy the same LDP, by means of the approximation lemma in \cite{Djellout1} and of three strong tools: Gronwall's inequality, L\'evy's inequality and an isoperimetric inequality for gaussian processes. By the estimation above (\ref{Drift1}) and (\ref{Drift3}), and as $Q_{\ell,t}^n(X-Y)$  was also estimated in the proof Theorem 1.3 
in \cite{Djellout1}. It remains to control $Z_{\ell,n}$ given in (\ref{defZ}) . For this we just need to prove that:
\begin{equation} 
\label{Drift_2}
\dfrac{1}{\varepsilon(n)}\dfrac{\sqrt{n}}{b_n}Z_{\ell,n}\superexpmdp 0.
\end{equation}

Chose $\varepsilon(n)>0$ so that 
\begin{equation}\label{epss}
\frac{\varepsilon(n)\sqrt{n}}{b_n} \to 0 \quad {\rm but} \quad \frac{\left(\frac{1}{n^2}+\eta\left(\frac{1}{n}\right)^2\right)b_n}{\varepsilon(n)\sqrt{n}} \to 0.
\end{equation}

By \eqref{Drift_4}  and the definition of $Z_{\ell,n}$ given in (\ref{defZ}), we have that
$$
\limsup_{n \to \infty}\dfrac{1}{b_n^2} \log \mathbb{P}\left(\dfrac{1}{\varepsilon(n)}\dfrac{\sqrt{n}}{b_n}Z_{\ell,n}>\delta \right)\leqslant \max(A,B) 
$$
where
\begin{align}\label{A_2}
A &=\limsup_{n \to \infty}\dfrac{1}{b_n^2} \log \mathbb{P}\left(\dfrac{1}{\varepsilon(n)b_n\sqrt{n}}\max_{k\leqslant n}\sup_{t_{k-1}^n \leqslant t \leqslant t_k^n}\|M_t-M_{t_{k-1}^n}\|^2>C_4\delta \right)\nonumber\\
&\leqslant2\max_{\ell=1,2} \limsup_{n \to \infty}\dfrac{1}{b_n^2} \log \mathbb{P}\left(\dfrac{1}{\varepsilon(n)b_n\sqrt{n}}\max_{k\leqslant n}\sup_{t_{k-1}^n \leqslant t \leqslant t_k^n}|M_{\ell,t}-M_{\ell,t_{k-1}^n}|^2>C_4\delta \right),
\end{align}
and
\begin{align}\label{B_2}
B &=\limsup_{n \to \infty}\dfrac{1}{b_n^2} \log \mathbb{P}\left(\dfrac{1}{\varepsilon(n)b_n\sqrt{n}}\left(\frac{1}{n^2}+\eta\left(\frac{1}{n}\right)^2\right)\sup_{0 \leqslant t \leqslant 1}\|M_t\|^2>C_5\delta \right)\nonumber\\
&\leqslant 2\max_{\ell=1,2}\limsup_{n \to \infty}\dfrac{1}{b_n^2} \log \mathbb{P}\left(\dfrac{1}{\varepsilon(n)b_n\sqrt{n}}\left(\frac{1}{n^2}+\eta\left(\frac{1}{n}\right)^2\right)\sup_{0 \leqslant t \leqslant 1}|M_{\ell,t}|^2>C_5\delta \right).
\end{align}

By L\'evy's inequality for a Brownian motion and our choice \eqref{epss} of $\varepsilon(n)$, the limit \eqref{B_2} are also $-\infty$.  As in \cite{Djellout1}, it's more little difficult to estimate \eqref{A_2}. By the isoperimetric inequality [\cite{Ledoux1}, p17,(1.24)] and our choice \eqref{epss},  we conclude that the limit  \eqref{A_2} are both $-\infty$.
\vspace{10pt}


\subsection{Proof of  Corollary \ref{LDPrhoconstant}}$\qquad$

\vspace{10pt}
We have just to do the identification of the rate function.   We knew  that $\varrho_1^n(X)$ satisfies the LDP on $\mathbb{R}$ with speed $n$ and with the good rate function given by
$$
 I^{\rho}_{ldp}(u) :=\inf_{\{(x_1,x_2,x_3)\in\mathbb{R}^3: x_3=u\sqrt{x_1x_2},x_1>0,x_2>0\}} I^{V}_{ldp}(x_1,x_2,x_3),$$
 where $I^{V}_{ldp}$ is given in (\ref{LDPtauxconstantV}). So
  $$
 I_{ldp}^{\rho}(u)=\inf\left\{ \log\left(\frac{\sqrt{\sigma_1^2\sigma_2^2(1-\rho^2)}}{\sqrt{x_1x_2}\sqrt{1-u^2}}\right)-1+\frac{\sigma_2^2 x_1+\sigma_1^2 x_2-2\rho\sigma_1\sigma_2 u\sqrt{x_1x_2}}{2\sigma_1^2\sigma_2^2(1-\rho^2)}, \quad x_1>0,x_2>0\right\}.
$$

 The above infinimum is attained at the point $(x_1,x_2)=\left(\dfrac{\sigma_1^2(1-\rho^2)}{1-\rho u},\dfrac{\sigma_2^2(1-\rho^2)}{1-\rho u}\right)$, so we obtain (\ref{LDPrhoconstanttaux1}).

 \vspace{10pt}


\subsection{Proof of Proposition \ref{MDP3}}$\qquad$

\vspace{10pt}
As said before, quite unusually, the MDP is here a little bit harder to prove, due to the fact that it is not a simple transformation of the MDP of $\frac{\sqrt n}{b_n}(V^n_t-[V]_t)$. Therefore we will use the strategy developped for the TCL: the delta-method. Fortunately, Gao and Zhao \cite{Gao1} have developped such a technology at the large deviations level. However it will require to prove quite heavy exponential negligibility to be able to do so. For simplicity we omit $X$ in the notations of  $Q_{1,t}^n(X)$ and $C_{t}^n(X)$.
\vspace{10pt}

Let introduce $\Xi_{t,n}$ such that $\Xi_{t}^n:=\sqrt{Q_{1,t}^n}\sqrt{Q_{2,t}^n}$. Then by the Lemma \ref{Gao} applied to the functions $g:=(x,y,z) \mapsto \sqrt{x}\sqrt{y}$ and  $h:=(x,y,z) \mapsto \frac{1}{\sqrt{x}\sqrt{y}}$, we deduce that $\frac{\sqrt{n}}{b_n}\biggl(\Xi_{1}^n-\mathbb{E}(\Xi_{1}^n)\biggr)$ and $\frac{\sqrt{n}}{b_n}\biggl((\Xi_{1}^n)^{-1}-({\mathbb{E}(\Xi_{1}^n)})^{-1}\biggr)$
 satisfies the LDP on $\mathbb{R}$ with the same speed $b_n^2$ and with the rates functions respectively given by $I_{mdp}^{\Xi}$ and $I_{mdp}^{\Xi^{-1}}$:
$$I_{mdp}^{\Xi}(u):=\inf_{\{(x,y,z)\in \mathbb{R}^3,u=\frac{\sigma_1^2y+ \sigma_2^2x}{2\sigma_1\sigma_2}\}}\{I_{mdp}(x,y,z)\},$$
and
$$I_{mdp}^{\Xi^{-1}}(u):=\inf_{\{(x,y,z)\in \mathbb{R}^3,u=-\frac{\sigma_1^2y+ \sigma_2^2x}{2\sigma_1^3\sigma_2^3}\}}\{I_{mdp}(x,y,z)\},$$
where $I_{mdp}$ is given in (\ref{MDPtaux1}).

By some simple calculations, we have
\begin{equation}\label{MDP31}\varrho_1^n(X)-\varrho =\aleph_1^n +\aleph_2^n+ \aleph_3^n+\aleph_4^n-\aleph_5^n -\aleph_6^n,
\end{equation}
where
$$\aleph_1^n:=(C_1^n-\mathbb{E}C_1^n)\biggl(\frac{1}{\Xi_{1}^n}-\frac{1}{\mathbb{E}(\Xi_{1}^n)}\biggr),\qquad\aleph_2^n:=(C_1^n-\mathbb{E}C_1^n)\frac{1}{\mathbb{E}\Xi_{1}^n},$$
$$\aleph_3^n:=(\mathbb{E}C_1^n-\varrho\mathbb{E}\Xi_{1}^n)\biggl(\frac{1}{\Xi_{1}^n}-\frac{1}{\mathbb{E}\Xi_{1}^n}\biggr),\qquad\aleph_4^n :=(\mathbb{E}C_1^n-\varrho\mathbb{E}\Xi_{1}^n) \frac{1}{\mathbb{E}\Xi_{1}^n},$$
$$\aleph_5^n :=\varrho(\Xi_{1}^n-\mathbb{E}\Xi_{1}^n)\biggl(\frac{1}{\Xi_{1}^n}-\frac{1}{\mathbb{E}\Xi_1^n}\biggr),\qquad\aleph_6^n :=\varrho(\Xi_{1}^n-\mathbb{E}\Xi_{1}^n)\frac{1}{\mathbb{E}\Xi_{1}^n}.$$

To prove the Theorem \ref{MDP3}, we have to use  the Lemma \ref{Gao} and prove some negligibility in the sence of MDP:  
\begin{equation}
\label{negli1}
\dfrac{\sqrt{n}}{b_n}\aleph_1^n\superexpmdp 0,
\end{equation}
\begin{equation}
\label{negli2}
\dfrac{\sqrt{n}}{b_n}\aleph_3^n\superexpmdp 0,
\end{equation}
\begin{equation}
\label{negli3}
\dfrac{\sqrt{n}}{b_n}\aleph_4^n\superexpmdp 0,
\end{equation}
\begin{equation}
\label{negli4}
\dfrac{\sqrt{n}}{b_n}\aleph_5^n\superexpmdp 0.
\end{equation}

Since $\mathbb{E}C_1^n-\varrho\mathbb{E}\Xi_{1,n} \to 0$ as $n \to \infty$, \eqref{negli3} follows.

We have for all $\delta>0$
$$\mathbb{P}\left(\frac{\sqrt{n}}{b_n}\left|\aleph_1^n\right|\geqslant \delta \right)\leqslant \mathbb{P}\left(\frac{\sqrt{n}}{b_n}\biggl|C_1^n-\mathbb{E}C_1^n\biggr| \geqslant \alpha_n\right) + \mathbb{P}\left(\frac{\sqrt{n}}{b_n}\biggl|\frac{1}{\Xi_{1,n}}-\frac{1}{\mathbb{E}\Xi_{1,n}}\biggr| \geqslant  \alpha_n \right),$$
where $\alpha_n=\sqrt{\frac{\sqrt{n}}{b_n}\delta}.$

So, by the Lemma 1.2.15 in \cite{Demzei}, we have that for all $\delta>0$
$$\limsup_{n \to \infty}\dfrac{1}{b_n^2} \log \mathbb{P}\left(\frac{\sqrt{n}}{b_n}\left|\aleph_1^n\right| \geqslant \delta\right)$$
is majorized by the maximum of the following two limits
\begin{equation}\label{C1}\limsup_{n \to \infty}\dfrac{1}{b_n^2} \log \mathbb{P}\left(\frac{\sqrt{n}}{b_n}\biggl|C_1^n-\mathbb{E}C_1^n\biggr| \geqslant \alpha_n \right),
\end{equation}
\begin{equation}\label{C2}\limsup_{n \to \infty}\dfrac{1}{b_n^2} \log \mathbb{P}\left(\frac{\sqrt{n}}{b_n}\biggl|\frac{1}{\Xi_{1,n}}-\frac{1}{\mathbb{E}\Xi_{1,n}}\biggr| \geqslant \alpha_n \right).
\end{equation}

Let $A>0$ be arbitrary, since $\alpha_n\rightarrow\infty$ as $n\to \infty$, so for $n$ large enough we obtain that 
$$\dfrac{1}{b_n^2} \log \mathbb{P}\left(\frac{\sqrt{n}}{b_n}\biggl|C_1^n-\mathbb{E}C_1^n\biggr| \geqslant \alpha_n\right) \leqslant \dfrac{1}{b_n^2} \log \mathbb{P}\left(\frac{\sqrt{n}}{b_n}\biggl|C_1^n-\mathbb{E}C_1^n\biggr| \geqslant A\right).$$

By the MDP of $\frac{\sqrt{n}}{b_n}(C_1^n-\mathbb{E}C_1^n)$ obtained in Theorem 2.3 in \cite{Djellout2}, and by letting $n$ to infinity, we obtain that 
$$\limsup_{n \to \infty}\dfrac{1}{b_n^2} \log \mathbb{P}\left(\frac{\sqrt{n}}{b_n}\biggl|C_1^n-\mathbb{E}C_1^n\biggr| \geqslant \alpha_n \right) \leqslant -\inf_{|x| \geqslant A} I_{mdp}^{C}(x).
$$

Letting $A$ gos to the infinity, we obtain that the term in (\ref{C1}) goes to $-\infty$.

By the MDP of $\frac{\sqrt{n}}{b_n}(\frac{1}{\Xi_{1,n}}-\frac{1}{\mathbb{E}(\Xi_{1,n})})$ stated before and in the same way we obtain that the term in (\ref{C2}) goes to $-\infty$. So we obtain (\ref{negli1}).

The same calculations give us \eqref{negli2} and \eqref{negli4}.

So 
$$\frac{\sqrt{n}}{b_n}(\varrho_1^n(X)-\varrho)$$
and
$$\frac{\sqrt{n}}{b_n}\biggl(C_1^n-\mathbb{E}C_{1}^n-\varrho(\Xi_{1}^n-\mathbb{E}\Xi_{1}^n)\biggl) \frac{1}{\mathbb{E}\Xi_{1}^n}$$
satisfies the same MDP. 

Since $\mathbb{E}(\Xi_{1,n})\longrightarrow \sigma_1\sigma_2$, so
$$\frac{\sqrt{n}}{b_n}(\varrho_1^n-\varrho)\quad {\rm and}\quad \frac{\sqrt{n}}{b_n}\biggl(C_1^n-\mathbb{E}C_{1}^n-\varrho(\Xi_{1}^n-\mathbb{E}\Xi_{1}^n) \biggr)\frac{1}{\sigma_1\sigma_2}
$$
satisfies the same MDP.
 
Then by the Lemma \ref{Gao} applied to the function  $\phi: (x,y,z) \mapsto (z-\varrho\sqrt{x}\sqrt{y})/\sigma_1\sigma_2$, we deduce that  $\frac{\sqrt{n}}{b_n}(\varrho_1^n(X)-\varrho)$ satisfies the LDP on $\mathbb{R}$ with speed $b_n^2$ and with the rate function given by
$$I^{\varrho}_{mdp}(u)=\inf_{\{(x,y,z)\in\mathbb{R}^3: u=\frac{z}{\sigma_1\sigma_2} - \varrho\frac{\sigma_1^2y+x\sigma_2^2}{2\sigma_1^2\sigma_2^2} \}}I_{mdp}(x,y,z),
$$
where $I_{mdp}$ is given in (\ref{MDPtaux1}).
\vspace{10pt}

\subsection{Proof of Proposition \ref{MDP4}}$\qquad$

\vspace{10pt}
By the Lemma \ref{Gao} applied to the function $f:x \mapsto \frac{1}{x}$, $\frac{\sqrt{n}}{b_n}(\frac{1}{Q_{1,1}^n}-\frac{1}{\mathbb{E}(Q_{1,1}^n)})$
 satisfies the LDP on $\mathbb{R}$ with the same speed $b_n^2$ and with the rate function given by 
 $$I_{mdp}^{Q_1^{-1}}(u):=\inf_{\{(x,y,z)\in\mathbb{R}^3: u=-\frac{x}{\sigma_1^4} \}}\{I_{mdp}(x,y,z)\}$$
 
By some simple calculations, we have
\begin{equation}\label{MDP41}
\beta_{1,1}^n(X)-\varrho\frac{\sigma_2}{\sigma_1} =\jmath_1^n+\jmath_2^n+\jmath_3^n+\jmath_4^n-\jmath_5^n -\jmath_6^n,
\end{equation}
where
$$\jmath_1^n:=\left(C_1^n-\mathbb{E}C_1^n\right)\biggl(\frac{1}{Q_{1,1}^n}-\frac{1}{\mathbb{E}Q_{1,1}^n}\biggr),\qquad \jmath_2^n:=(C_1^n-\mathbb{E}C_1^n)\frac{1}{\mathbb{E}Q_{1,1}^n},
$$
$$\jmath_3^n:=(\mathbb{E}C_1^n-\varrho\frac{\sigma_2}{\sigma_1}\mathbb{E}Q_{1,1}^n)\biggl(\frac{1}{Q_{1,1}^n}-\frac{1}{\mathbb{E}Q_{1,1}^n}\biggr),\qquad \jmath_4^n:=(\mathbb{E}C_1^n-\varrho\frac{\sigma_2}{\sigma_1}\mathbb{E}Q_{1,1}^n)\frac{1}{\mathbb{E}Q_{1,1}^n},
$$
$$\jmath_5^n:=\varrho\frac{\sigma_2}{\sigma_1}\left(Q_{1,1}^n-\mathbb{E}Q_{1,1}^n\right)\biggl(\frac{1}{Q_{1,1}^n}-\frac{1}{\mathbb{E}Q_{1,1}^n}\biggr),\qquad\jmath_6^n:=\varrho\frac{\sigma_2}{\sigma_1}\left(Q_{1,1}^n-\mathbb{E}Q_{1,1}^n\right)\frac{1}{\mathbb{E}Q_{1,1}^n}.$$

To prove the Theorem \ref{MDP4}, we have to use  the Lemma \ref{Gao} and prove some negligibility in the sense of MDP:

\begin{equation}
\label{negli5}
\dfrac{\sqrt{n}}{b_n}\jmath_1^n\superexpmdp 0,\quad\dfrac{\sqrt{n}}{b_n}\jmath_3^n\superexpmdp 0,\quad \dfrac{\sqrt{n}}{b_n}\jmath_4^n\superexpmdp 0,\quad
\dfrac{\sqrt{n}}{b_n}\jmath_5^n\superexpmdp 0.
\end{equation}

The same calculations as for the negligibility of $\aleph_j^n$ works here to obtain (\ref{negli5}).

Since $\mathbb{E}Q_{1,1}^n\longrightarrow \sigma_1^2$, we deduce that
$$\dfrac{\sqrt{n}}{b_n}(\beta_{1,1}^n(X)-\varrho\frac{\sigma_2}{\sigma_1}))$$
and
$$\dfrac{\sqrt{n}}{b_n}(\left(C_1^n-\mathbb{E}C_1^n- \varrho\frac{\sigma_2}{\sigma_1}(Q_{1,1}^n-\mathbb{E}Q_{1,1}^n\right)\frac{1}{\sigma_1^2}$$
satisfies the same MDP.
 
Then by the Lemma \ref{Gao} applied to the function $k: (x,y,z) \mapsto (z-\varrho\frac{\sigma_2}{\sigma_1}x)/\sigma_1^2$ we deduce that  $\frac{\sqrt{n}}{b_n}(\beta_{1,1}^n(X)-\varrho\frac{\sigma_2}{\sigma_1})$ satisfies the LDP on $\mathbb{R}$ with speed $b_n^2$ and with the rate function given by
 $$I_{mdp}^{\beta_{1,1}}(u)=\inf_{\{(x,y,z)\in\mathbb{R}^3: u=(z- \varrho\frac{\sigma_2}{\sigma_1}x)/\sigma_1^2 \}}I_{mdp}(x,y,z)$$

where $I_{mdp}$ is given in (\ref{MDPtaux1}).


\section{Appendix}
The proofs of the LDP in Theorems \ref{LDP} and \ref{LDP2} are respectively based on the Lemmas \ref{Najim_Jamal_1} and \ref{Najim12} that we will present here for completeness.\\\vspace{10pt}

\subsection{Avoiding G\"artner-Ellis theorem by Najim \cite{Najim1,Najim3}}

Let us introduce some notations and assumptions in this section.

Let $\mathcal{X}$ be a topological vector compact space endowed with it's Borel $\sigma-$field $\mathcal{B}(\mathcal{X})$. Let $\mathcal{B}V([0,1],\mathbb{R}^d),$ (shortened in $\mathcal{B}V$)
be a space of functions of bounded variation on $[0,1]$ endowed with it's Borel $\sigma-$field $\mathcal{B}_w$. Let ${\mathcal P}(\mathcal X)$ the set of probability measures on $\mathcal X$.

Let $\tau(z)=e^{|z|}-1, z\in \mathbb{R}^d$ and let us consider
\begin{eqnarray*}
 \mathcal{P}_{\tau}(\mathbb{R}^d) &= &\left\{P\in \mathcal{P}(\mathbb{R}^d), \exists a>0; \int_{\mathbb{R}^d}\tau\left(\frac{z}{a}\right) P(dz) < \infty\right\}\\
&=& \left\{P\in \mathcal{P}(\mathbb{R}^d), \exists \alpha>0; \int_{\mathbb{R}^d}e^{\alpha|z|} P(dz) < \infty\right\}.
\end{eqnarray*}

$\mathcal{P}_{\tau}$ is the set of probability distributions having some exponential moments.We denote by $M(P,Q)$ the set of all laws on $\mathbb{R}^d\times\mathbb{R}^d$ with given
marginals $P$ and $Q$. We introduce the Orlicz-Wasserstein distance defined on $\mathcal{P}_{\tau}(\mathbb{R}^d)$ by
$$d_{OW}(P,Q)= \inf_{\eta \in M(P,Q)} \inf \biggl\{a>0; \int_{\mathbb{R}^d\times\mathbb{R}^d}\tau\left(\frac{z-z'}{a}\right) \mathrm \eta(dzdz') \leqslant 1 \biggr\}
$$

Let $(Z_i^n)_{1 \leqslant i \leqslant n, n \in \mathbb{N}}$ be a sequence of $\mathbb{R}^d-$valued independent random variables satisfying:
\begin{itemize}
 \item [N-1]$$\mathbb{E}e^{\alpha.||Z||}< +\infty \qquad for\quad some \quad \alpha > 0.$$
 \item [N-2]Let $\left(x_i^n, 1 \leqslant i \leqslant n, n \geqslant 1\right)$ be a $\mathcal{X}-$valued sequence of elements satisfying:
    $$\dfrac{1}{n}\sum_{j=1}^{n}\delta_{x_j^n} \xrightarrow[{n \to \infty}]{weakly}\quad R.$$
     Where $R$ is Assumed to be a strictly positive probability measure, that is $R(U) > 0$ whenever $U$ is a nonempty open subset of $\mathcal{X}.$\\
 \item [N-3]$\mathcal{X}$ is a compact space.\\
 \item [N-4]There exist a family of probability measure $(P(x,\cdot),x\in \mathcal{X})$ over $\mathbb{R}^d$ and a sequence $\left(x_i^n, 1 \leqslant i \leqslant n, n \geqslant 1\right)$ with values in $\mathcal{X}$ such that 
       the law of each $Z_i^n$ is given by:
       $$\mathcal{L}(Z_i^n) \sim P(x_i^n,dz).$$
       We will equally write $P(x,\cdot), P_x \quad or\quad P_x(dz).$\\
 \item [N-5]Let $(P(x,\cdot),x\in \mathcal{X}) \subset \mathcal{P}_{\tau}(\mathbb{R}^d)$ be a given distribution. The application $x\mapsto P(x,A)$ is measurable
   whenever the set $A\subset \mathbb{R}^d$ is borel. Morever, the function
   \begin{align*}
    \Gamma:& \mathcal{X} \to \mathcal{P}_{\tau}(\mathbb{R}^d)\\
           & x\mapsto P(x,\cdot)
   \end{align*}

is continuous when $\mathcal{P}_{\tau}(\mathbb{R}^d)$ is endowed with the topology induced by the distance $d_{OW}$
\end{itemize}

\begin{lem} Theorem 2.2 in \cite{Najim1}
\label{Najim_Jamal_1}

Assume that $Z_i^n$ are independent and identically distributed, so we denote $Z_i^n$ by $Z_i$.

 Assume that (N-1) and (N-2) hold. Let $f:\mathcal{X}\to \mathbb{R}^{m\times d}$ be a (matrix-valued) bounded continuous function, such that
$$
f(x)\cdot z = \begin{pmatrix}
                                       \vspace{0.75cm}
                                       f_1(x)\cdot z\\
                                        \vspace{0.75cm}
                                       \vdots\\
                                       f_m(x)\cdot z\\
                                       \end{pmatrix}
$$
where each $f_j\in C_d(\mathcal{X})$ is the $j^{th}$ row of the matrix $f$.

Then the family of the weighted empirical mean 
$$\left\langle L_n,f \right\rangle: = \dfrac{1}{n}\sum_{i=1}^{n}f(x_i^n)\cdot Z_i$$
satisfies the LDP in $(\mathbb{R}^m,\mathcal{B}(\mathbb{R}^m))$ with speed $n$ and the good rate function
$$I_f(x)=\sup_{\theta\in\mathbb{R}^m}\{\left\langle \theta,x \right\rangle - \int_{\mathcal{X}}\Lambda[\sum_{i=1}^{m}\theta_i\cdot f_i(x)] \mathrm R(dx)\} \qquad \forall x\in \mathbb{R}^m$$
where $\Lambda$ denote the cumulant generating function of Z
$$\Lambda(\lambda) = \log\mathbb{E}e^{\lambda \cdot Z} \qquad for \quad \lambda \in \mathbb{R}^d$$
\end{lem}
\vspace{10pt}

\begin{lem} Theorem 4.3 in \cite{Najim2}
 \label{Najim12}
 
 Assume that (N-1), (N-2), (N-3), (N-4) and (N-5) hold. Then the family of random functions
 $$\overline Z_n(t)=\frac{1}{n}\sum_{k=1}^{[nt]}Z_k^n, \quad t\in [0,1]$$
 satisfies the LDP in $(\mathcal{BV},\mathcal{B}_w)$ with the good rate function
 $$\phi(f)=\int_{[0,1]} \Lambda^*(x,f_a^{'}(x)) \mathrm dx + \int_{[0,1]}\rho(x,f_s^{'}(x)) \mathrm d\theta(x)$$
 where $\theta$ is any real-valued nonnegative measure with respect to which $\mu_s^f$ is absolutely continuous and $f_s^{'}=\frac{d\mu_s^f}{d\theta},$ where $$\Lambda^*(x,z)= \sup_{\lambda \in \mathbb{R}^d} \biggl\{\left\langle \lambda,z \right\rangle - \Lambda(x,\lambda) \biggr\}, \quad \forall z\in \mathbb{R}^d$$
 with $\Lambda(x,\lambda)=\log \int_{\mathbb{R}^d} e^{\left\langle \lambda,z \right\rangle}  P(x,dz), \quad \forall \lambda \in \mathbb{R}^d$ and the recession function $\rho(x,z)$ of $\Lambda^*(x,z)$ defined by: $\rho(x,z)=\sup\{\left\langle \lambda,z \right\rangle, \lambda \in D_x\}$ with $D_x=\{\lambda \in \mathbb{R}^d, \Lambda(x,\lambda)<\infty\}.$
\end{lem}
\vspace{10pt}

\subsection{Delta method for large deviations \cite{Gao1}} 

In this section, we recall the delta method in large deviation.

Let $\mathcal{X}$ and $\mathcal{Y}$ be two metrizable topological linear spaces. A function $\phi$ defined on a subset $\mathcal{D}_{\phi}$ of $\mathcal{X}$
with values on $\mathcal{Y}$ is called Hadamard differentiable at $x$ if there exists a continuous functions $\phi': \mathcal{X} \mapsto \mathcal{Y}$ such that
\begin{equation}
 \lim_{n\to \infty} \dfrac{\phi(x+t_nh_n)-\phi(x)}{t_n}=\phi'(h)
\end{equation}
holds for all $t_n$ converging to $0+$ and $h_n$ converging to h in $\mathcal{X}$  such that $x+t_nh_n \in \mathcal{D}_{\phi}$ for every $n$.

\begin{lem}
 \label{Gao}
 Let $\mathcal{X}$ and $\mathcal{Y}$ be two metrizable topological linear spaces. Let $\phi:\mathcal{D} \subset \mathcal{X} \mapsto \mathcal{Y}$ be a
 Hadamard differentiable at $\theta$ tangentially to $\mathcal{D}_0$, where $\mathcal{D}_{\phi}$ and $\mathcal{D}_0$ are two subset of $\mathcal{X}$. Let
  $\{(\Omega_n, \mathcal{F}_n, \mathbb{P}_n), n\geqslant 1\}$ be a sequence of probability space and let $\{X_n, n\geqslant 1\}$ be a sequence of maps from from $\Omega_n$ to $\mathcal{D}_{\phi}$ 
  and let $\{r_n, n\geqslant 1\}$ be a sequence of positive real numbers satisfying $r_n \to +\infty$ and let $\{\lambda(n), n\geqslant 1\}$ be a sequence of positive real numbers satisfying $\lambda(n) \to +\infty$.\\
  If $\{r_n(X_n-\theta), n\geqslant 1\}$ satifies the LDP with speed $\lambda(n)$ and rate function I and $\{I<\infty\} \subset \mathcal{D}_0$, then
  $\{r_n(\phi(X_n)-\phi(\theta)), n\geqslant 1\}$ satifies the LDP with speed $\lambda(n)$ and rate function $I_{\phi'_{\theta}}$, where
  \begin{equation}
   I_{\phi'_{\theta}}(y)=inf\{I(x);\phi'_{\theta}(x)=y\}, \qquad y\in \mathcal{Y}
  \end{equation}

\end{lem}

\medskip

\nocite{*}

\bibliographystyle{acm}
\bibliography{biblio}

\nocite{*}

\vspace{10pt}

\end{document}